\documentclass[journal]{IEEEtran}

\usepackage{xcolor,soul,framed} 

\colorlet{shadecolor}{yellow}
\usepackage[pdftex]{graphicx}
\graphicspath{{../pdf/}{../jpeg/}}
\DeclareGraphicsExtensions{.pdf,.jpeg,.png}

\usepackage[cmex10]{amsmath}
\usepackage{array}
\usepackage{mdwmath}
\usepackage{mdwtab}
\usepackage{eqparbox}
\usepackage{url}
\usepackage{xcolor}

\usepackage{multirow}
\usepackage{algorithm,algpseudocode}
\usepackage{comment}
\usepackage{cancel}

\usepackage{multicol}
\usepackage{amssymb}
\usepackage{mathrsfs}
\usepackage{amsmath}
\usepackage[english]{babel}
\usepackage{float}
\usepackage{accents}
\usepackage{bm}
\usepackage{mathtools}
\usepackage{verbatim}
\usepackage{easy-todo}
\usepackage[standard]{ntheorem}

\newcommand{\bb}[1]{\mathbf{#1}}
\newcommand*{\QEDA}{\null\nobreak\hfill\ensuremath{\square}}%
\newcommand{\x}{\bb{x}}

\newcommand{\q}{\bb{q}}

\newcommand{\p}{\bb{p}}
\newcommand{\ut}{\bb{u}}
\newcommand{\R}{\mathbb{R}}

\newcommand{\nor}{\bb{n}}

\newcommand\mbf[1]{\mathbf{#1}}
\newcommand\norm[1]{\left\lVert#1\right\rVert}

\newcommand{\Ht}{L^{2}(0,T;H^1(\Omega))}
\newcommand{\Htd}{L^{2}(0,T;H^1(\Omega)^*)}

\algnewcommand{\And}{\textbf{and}}

\begin{document}
\bstctlcite{IEEEexample:BSTcontrol}
    \title{Density control of large-scale particles swarm through PDE-constrained optimization}
  \author{Carlo Sinigaglia,
    Andrea Manzoni,
    Francesco Braghin,~\IEEEmembership{Fellow,~IEEE}

  \thanks{Carlo Sinigaglia is a PhD Candidate at Politecnico di Milano, Department of Mechanical Engineering, Milano 20133, Italy (e-mail: carlo.sinigaglia@polimi.it), Corresponding author.} 
  \thanks{Prof. Andrea Manzoni is Associate Professor at Politecnico di Milano, MOX - Department of Mathematics, Milano 20133, Italy (e-mail: andrea1.manzoni@polimi.it).}
  \thanks{Prof. Francesco Braghin is Full Professor at Politecnico di Milano, Department of Mechanical Engineering, Milano 20133, Italy (e-mail: francesco.braghin@polimi.it). }
  }


\maketitle

\begin{abstract}

We describe in this paper an optimal control strategy for shaping a large-scale swarm of particles using boundary global actuation. This problem arises as a key challenge in many swarm robotics applications, especially when the robots are passive particles that need to be guided by external control fields. The system is large-scale and underactuated, making the control strategy at the microscopic particle level infeasible. We consider the Kolmogorov forward equation associated to the stochastic process of the single particle to encode the macroscopic behaviour of the particles swarm. The control inputs shape the velocity field of the density dynamics according to the physical model of the actuators. We find the optimal actuation considering an optimal control problem whose state dynamics is governed by a linear parabolic advection-diffusion equation where the control induces a transport field.
From a theoretical standpoint, we show the existence of a solution to the resulting nonlinear optimal control problem.
From a numerical standpoint, we employ the discrete adjoint method to accurately compute the reduced gradient and we show how it commutes with the optimize-then-discretize approach. Finally, numerical simulations show the effectiveness of the control strategy in driving the density sufficiently close to the target.
\end{abstract}


%
\IEEEpeerreviewmaketitle


\section{Introduction}
\label{intro}

\IEEEPARstart{L}{arge-scale} underactuated robotic systems in the form of particle swarms are increasingly finding applications in robotics. Micro-robots driven by the uniform magnetic field generated by a Magnetic Resonance Imaging (MRI) system are foreseen to be used in medical applications \cite{boundary_1,boundary_2} to control swarms of robots for drug delivery in the human body. 

Particles swarms find also application in distributed space robotics. A large-scale system of reflective particles is envisioned to replace continuum monolithic apertures for advancing the current state of the art space telescope technology \cite{NIAC,DCD,granular_dyn} where a cloud of reflective particles is deployed and kept in shape by external electromagnetic actuators; the conceptual design, feasibility analysis and working principles have been studied in \cite{optics_gran,granular_dyn,Unc_gran}. A fundamental advance in the field of micro-robotics appeared very recently in \cite{Nature2020} where a remote laser-controlled actuation system has been developed together with a method to mass produce a swarm of robots in the order of millions of agents. 

In this paper, we provide a method to control such large-scale and underactuated systems by steering the macroscopic density dynamics of the particles using methods from optimal control theory, where the state dynamics is described by a partial differential equation. This framework allows us to find the boundary actuation time history that steers the swarm towards a given target distribution. The velocity vector field in the state dynamics is determined by boundary controls according to the actuator physical model that ends up in a linear combination of control functions. As a consequence, the control functions enter bilinearly in the state equation, making the Optimal Control Problem (OCP) nonlinear.

The robotic swarms that we consider cannot perform any onboard computation and moves according to external stimuli. We model the robot as a particle subjected to a velocity field that can be externally modulated as a function of the boundary actuation. Therefore, all the particles are subjected to the same control field. The control authority is limited to the boundary of the workspace and cannot be modulated arbitrarily in space and time to guide and control each single particle.
From the control-theoretic perspective, this results in a largely underactuated system where we aim at controlling the macroscopic density instead of following the dynamics of each single particle.

\subsection{Literature survey on control of swarms of particles}

The existing literature on density control of large-scale multi-agent systems is mostly based on robotics swarm, where each agent is able to directly exert actions to control its motion. Since the control problem becomes intractable as the number of individual robots gets large, control laws based on macroscopic models of such systems have recently become an attractive alternative to classical optimal control and path planning methods \cite{karthik_survey}.
The main idea of these methods is to optimize the macroscopic behavior based on a density evolution model that takes the form of an unsteady Advection-Diffusion (AD) equation. 
In \cite{Milutinovic2006ModelingAO} Pontryagin's Minimum Principle is used to recover optimal density motion and task allocation in a stochastic hybrid automation one dimensional setting.

In \cite{gen_grad,nc_fod}, the necessary optimality conditions were established using methods from Calculus of Variations, and then numerically solved for the case in which the control functions are null at the boundary. In \cite{kart_1,kart_2}, the same problem is tackled from the perspective of functional analysis stating some well-posedness and existence results. The controllability of the Advection-Diffusion equation using the velocity field as control input was addressed in \cite{karthik_bil} where finite-time and path controllability are proved for unconstrained velocity control fields. A similar PDE macroscopic model is analysed in \cite{rob_dep} and general well-posedness and existence results are obtained adapting a classical proof from \cite{lieberman}. Besides robotic applications, the problem of controlling the advection field for a steady Advection-Diffusion equation is also considered in \cite{Ito_97} where a regularized least-square identification problem is solved for an unknown divergent-free velocity field. A one-dimensional OCP in the advection field is considered in \cite{joshi_05} for an unsteady Advection-Diffusion equation with Dirichlet boundary conditions, for which the control field is allowed to vary both in space and time. 
A divergence free, space-time dependent control vector field is tackled more recently in \cite{glowinski2021bilinear} where an existence theorem for $L^2$ controls is provided in the case of Dirichlet boundary conditions, then an efficient conjugate gradient method is also developed for the numerical solution of the problem.

All the previously mentioned approaches are indeed attractive for robotic systems where the local controllers can follow without constraints on the optimal velocity vector field. However, in our case, the velocity vector field is constrained before-hand and is generated by boundary actuation, making the previous methods not suitable. Boundary global actuation is considered in \cite{boundary_1,boundary_2} for positioning a system of two particles using magnetic actuation. The particles are confined into a two-dimensional squared workspace and are subject to a uniform global control input, each particle experiences the same displacement and the symmetry is broken by assuming that particles pushed against the boundary remains still. However, this method does not scale well when the number of particles increases. 

On the other hand, the reformulation as a density control problem allows us to remove such a boundary constraint and to develop a control algorithm that is independent of the number of particles. In this paper, we will just assume that the particles cannot leave the domain considering no-flux boundary conditions. A numerical optimization algorithm is then able to find the actuation history that drives the density towards the target.

Compared to control strategies that make use of a PDE model, in our case the space dependence of the velocity vector field does not arise as an approximation of the control field as in \cite{gen_grad,nc_fod}, but it is imposed by the physical model of the actuators. 
In particular, the velocity field $\bb{v}$ depends on a set of control functions that allow us to ensure a nonnegativity constraint stemming from the physical nature of the actuators. These latter can only provide unilateral actions, thus limiting the control authority on the induced vector field.

Finally, we note that an existence result for a similar OCP is shown in \cite{Fleig2017}, where the velocity field is affine in the controls, however, for a Fokker-Planck equation involving Dirichlet boundary conditionswhile Neumann no-flux boundary conditions are considered in \cite{roy2018fokker} where, however, the control vector field nullify at the boundaries and the cost functional is linear in the state variable thus simplifying the form of the adjoint equations. Building on this work, we extend the existence result considering a Kolmogorov equation equipped with Neumann boundary conditions where the velocity vector field is induced by control functions defined at the boundary of the workspace.

\subsection{Main contributions and paper organization}

The main contributions of this paper, compared to the existing literature, are threefold:
\begin{enumerate}
  \item a suitable PDE model of the density evolution that is consistent with the particles dynamics is analysed; 
  
  \item an OCP to steer the probabilistic density of the particles towards a target distribution is set up. The boundary control field determines the velocity that drives the state dynamics. The existence of an optimal control in space and time is proven under box constraints imposed by the actuation model, then a set of first-order necessary optimality conditions is derived. Afterwards, a numerical approximation with respect to both space and time leads to a discrete set of state and adjoint equations;

  \item the discrete adjoint method is then used to efficiently compute the reduced gradient, and is shown how it commutes with its discretized continuous counterpart. Numerical simulations are finally performed to show the effectiveness of the algorithm considering different target functions of increasing difficulty.
\end{enumerate}

The paper is organized as follows. In Section~\ref{density_dyn}, the density dynamics is established together with its dependence on the control functions. In Section~\ref{OCP}, the OCP is set up, and an existence theorem is proved. Then, a set of first-order optimality conditions is derived. In Section~\ref{Num_approx}, numerical approximations of the system of optimality conditions is carried out in view of the numerical solution of the OCP. Then, the discrete adjoint method is compared to the ``Optimize then Discretize" approach. Numerical results are then shown in Section~\ref{num_sim}, while some conclusions then follow in Section~\ref{conclusion}.

\section{Density dynamics: the state problem}
\label{density_dyn}
In this section we establish the single particle stochastic model providing in detail the actuator layouts and their effect on the motion of the particles. Then, the density dynamics in the form of a time-dependent Advection-Diffusion partial differential equation is derived, and some properties of the solution of this equation are proved.

\begin{figure}
 \centering
 \includegraphics[width = \linewidth]{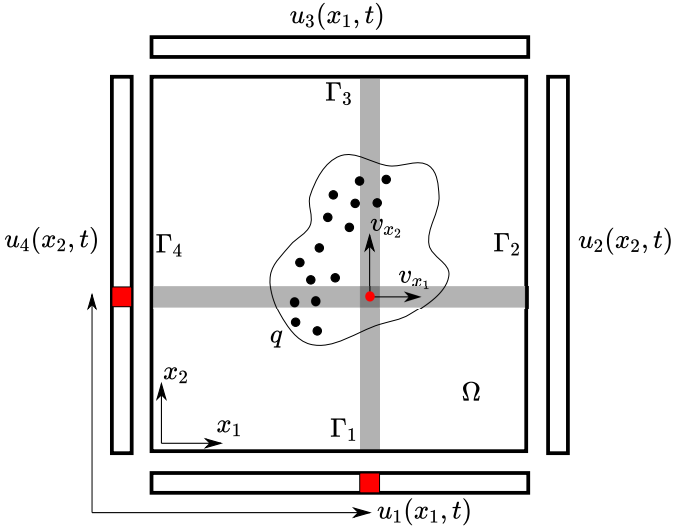}
 \caption{Layout of the control system, the actuators on the four sides generate a velocity field that determines the evolution of the density inside the workspace. The control functions are the intensities of the actuators. The boundary actuation induces a velocity field on all the particles located in the respective row or column.}
 \label{stoch_part}
 \end{figure}

We consider a simplified two-dimensional layout given, for the sake of simplicity, by the unit square, where the swarm of particles is confined as shown in Figure \ref{stoch_part}. The actuator stacks are positioned at each of the four sides and exert their action in the orthogonal direction with vanishing intensity as the distance from the relative source grows. The same layout was considered in \cite{OT} and in \cite{chinese_ot}, although in these works the modeling was carried out at the discrete level from the beginning.
Under the influence of the force field generated by the actuators, we assume that the particles motion follows the stochastic differential equation (SDE) of a reflective diffusion process \cite{karthik_bil}. 

\noindent Consider a swarm of $N$ particles that are deployed on the domain $\Omega=(0,1)^2$ with boundary $\Gamma = \partial \Omega$. The position of each particle is denoted by $\bb{X}_i(t)$, where $t$ denotes time and $i$ the particle's index. We assume that the position of each particle evolves according to a stochastic process and that the particles are noninteracting. Therefore, the random variables corresponding to the particle dynamics are independent and identically distributed, hence, we can drop the subscript $i$ and consider a single stochastic process $\bb{X}(t) \in \Omega$. The boundary actuation induces a deterministic velocity field $\bb{v}(\bb{X}(t),t) \in \R^2$, while the motion is perturbed by a two-dimensional Wiener process $\bb{W}(t)$.

This model is able to capture simultaneously the stochastic effect due to Gaussian diffusion, the deterministic motion imposed by the actuators and the reflection due to the boundedness of the workspace. The equation governing the microscopic particle dynamics is:
\begin{equation}
\label{micro_dyn}
\begin{cases}
d \mathbf{X}(t) &=\mathbf{v}(\mathbf{X}(t), t) d t+\sqrt{2 \mu} d \mathbf{W}(t)+\mathbf{n}(\mathbf{X}(t)) d \psi(t) \\
\mathbf{X}(0) &=\mathbf{X}_{0} ,
\end{cases}
\end{equation}
where $\bb{X}_{0} \in \R^2$ is the initial position of the particle, $\psi(t) \in \R$ is the reflecting function or local time (see e.g. \cite{karthik_bil}) that constrains $\bb{X}(t)$ to belong to the domain $\Omega$, $\mu>0$ is the diffusion constant and $\bb{n}(\bb{x})$ is the outward normal at $\bb{x}\in \partial \Omega$.
Equation (\ref{micro_dyn}) is the same model used in \cite{kart_1,kart_2,karthik_bil,rob_dep} to describe the microscopic dynamics of a robotic agent.
Considering the swarm of particles as a continuum, the density dynamics is described by the following Kolmogorov forward equation:
\begin{equation}
\label{state_eq}
\begin{aligned}
 \quad &   \frac{\partial q}{\partial t} +  \nabla \cdot (-\mu\nabla q + \bb{v} q)   =0   \hspace{0.95cm}  in \,\, \Omega \times (0,T), \\
& (-\mu \nabla q + \bb{v} q)\cdot \nor = 0 \hspace{1.9cm} on \,\, \Gamma \times (0,T),
\phantom{space}  \\
&                         q(\x,0) = q_0(\x)  \hspace{2.8cm} in \,\, \Omega \,\, at \,\, t=0 \\
\end{aligned}
\end{equation}
where $q \colon \Omega \times [0,T] \mapsto \R$ represents the probability density of the particles. Equation (\ref{state_eq}) is related to the SDE (\ref{micro_dyn}) through the relation $\mathbb{P}(\bb{X}(t) \in A) = \int_{A} q(\x,t) \, d\Omega$ for all $t \in [0,T]$ and all measurable $A \subset \Omega$. Note that Neumann conditions are imposed on the whole boundary, prescribing a no-flux condition. 

\begin{remark}
Note that the density dynamics is consistent with the microscopic particle dynamics in the sense that the total mass is conserved, and the particles are not allowed to leave $\Omega$. Indeed, defining $M(t)=\int_{\Omega} \, q \, d \Omega$ as the total mass in the workspace, we have that:
\begin{equation*}
\begin{aligned}
\dot{M}(t) & =  \int_{\Omega}  \frac{\partial q}{\partial t} \, d\Omega = \int_{\Omega}  -\nabla \cdot (-\mu\nabla q + \bb{v} q) \, d\Omega  \\&= \int_{\Gamma} -(-\mu \nabla q + \bb{v} q)\cdot \nor \, d \Gamma = 0,
\end{aligned}
\end{equation*}
using the Divergence Theorem and the no-flux boundary condition. The state dynamics in (2) is a particular case of the more general conormal derivative problem for linear parabolic PDEs, see, e.g. \cite{lieberman,rob_dep,evans}. 
\end{remark}

\noindent The boundary control field induces a velocity field in the domain $\Omega$. We assume that the actuators have a vanishing effect $\propto e^{-cx}$ where $x$ stands for the distance from the actuator and $c \in \R$ is a constant scalar that defines the decay rate. Such modeling choice is an abstraction of the typical behaviour of magnetic actuators \cite{act_model}. However, we stress the fact that the functional form of the decay rate does not affect the structure of the problem, and different kind of actuating models (e.g. laser-based) can be used. Hence, the velocity field induced by the control functions takes the form:

\begin{equation}
\label{act}
\mathbf{v}(\x,t) = 
\begin{bmatrix}
u_4(x_2,t) e^{-cx_1} - u_2(x_2,t) e^{-c(1-x_1)} \\
u_1(x_1,t) e^{-cx_2} - u_3(x_1,t) e^{-c(1-x_2)} \\
\end{bmatrix}.
\end{equation}

\noindent  The components $ v_{1}(x_1,x_2,t)$ and $ v_{2}(x_1,x_2,t)$ of the velocity field are linear combinations of the actuator functions defined on the boundary, weighted by negative exponentials. We note that the vector field $\bb{v}$ inherits the regularity of the control functions $u_1,\ldots,u_4$.

\noindent The well-posedness of the state problem can be proved by adapting classical results found, e.g. in  \cite{lieberman,quarteroni_valli,evans}.

We will now describe the functional setting for our problem and prove some estimates that will be used in Section \ref{OCP} to show the existence of optimal controls.
We define the space of controls for each boundary control function as $\mathcal{U}_{ad,i}:=\left\{u \in \mathcal{U}_i \quad \text { s.t. } \quad 0 \leq u_{i}(t) \leq u_{\max } \quad \text { a.e. } \quad\right. \text { on } \quad \Gamma_i \}$ with $\mathcal{U}_i = L^2(0,T;L^{\infty}(\Gamma_i))$; then we group the boundary control functions in a vector $\mathbf{u}=(u_1\,,\,\ldots\,,\,u_4)$ with associated space $\mathcal{U}_{ad} = \mathcal{U}_{ad,1} \times \mathcal{U}_{ad,2} \times \mathcal{U}_{ad,3} \times \mathcal{U}_{ad,4}$ endowed with the norm $\norm{\mathbf{u}}_{\mathcal{U}} = \sqrt{ \sum_{i=1}^4 \norm{u_i}_{\mathcal{U}_i}^2}$. First of all we show that the norm of the velocity field is bounded by the norm of the control field.

\begin{lemma}[Estimate on the velocity field]
\label{lemma_1}
Assume that $\mathbf{u} \in \mathcal{U}_{ad} $. Then, the following inequality holds:
\begin{equation}
\norm{\bb{v}}^2_{L^2(0,T;L^{\infty}(\Omega)^2)} \leq  8\, \norm{\mathbf{u}}^2_{\mathcal{U}}.
\end{equation}
\qquad Proof: see the Appendix.
\QEDA
\end{lemma}
In the remainder of the paper we will rely on the weak formulation of the state problem (\ref{state_eq}). This is obtained by multiplying (\ref{state_eq}) by a test function $\phi \in H^{1}(\Omega)$ and integrating over $\Omega$. The Divergence Theorem is then applied to handle the boundary condition so that, for every $t>0$, we obtain the following problem: find $q \in H^{1}(0,T;H^1(\Omega),H^1(\Omega)^{*})$ such that $\forall \phi \in H^{1}(\Omega), \,\, a.e. \,\, t \in (0,T)$ it holds that
\begin{equation} 
\label{weak_form}
\begin{cases}
&  \int_{\Omega} \frac{\partial q(t)}{\partial t} \phi \, d \Omega +\int_{\Omega}(\mu \nabla q(t) \cdot \nabla \phi-\mathbf{v}(t) \cdot \nabla \phi \, q(t)) d \Omega= 0 \\
& q(0)=q_0.
\end{cases}
\end{equation}
We can now define the bilinear forms associated to this problem as:
\begin{equation}
\label{bil_form}
    a(q,\phi) = d(q,\phi) + b(q,\phi ; \bb{v}(t)) = \int_{\Omega}(\mu \nabla q \cdot \nabla \phi-\mathbf{v}(t) \cdot \nabla \phi \, q) \, d \Omega
\end{equation}
where $d(q,\phi) :=\int_{\Omega}\mu \nabla q \cdot \nabla \phi \, d \Omega$ is associated to the diffusion term, while $b(q,\phi ; \bb{v}(t)) :=- \int_{\Omega} \mathbf{v}(t) \cdot \nabla \phi \, q \, d \Omega $ is associated to the advection term. The details on the derivation of the bilinear forms \eqref{bil_form} is given in the Appendix.

\begin{remark}
Note that the weak formulation can be expressed equivalently as: find $q$ such that $\forall \phi \in  H^{1}(\Omega), \,\, a.e. \,\, t \in (0,T)$
\begin{equation}
\begin{aligned}
\label{state_weak_2}
& \int_{\Omega} \frac{\partial q(t)}{\partial t} \phi + \mu \nabla q(t) \cdot \nabla \phi + \bb{v}(t)\cdot \nabla q \, \phi \, d\Omega \\
& +\int_{\Omega} (\nabla \cdot \bb{v}(t)) \, q(t) \, \phi \,\, d\Omega -\int_{\Gamma} \bb{v}(t) q(t) \cdot \bb{n} \, \phi \, d\Gamma = 0 \quad  \\
\end{aligned}
\end{equation}
and in this form it will be used for the numerical approximation in Section \ref{Num_approx}.
\end{remark}

The bilinear form $a(q,\phi)$ is weakly coercive, according to the following Lemma.
\begin{lemma}[Weak coercivity of $a(q,\phi)$]
\label{weakly_coerc}
Under the assumptions of Lemma \ref{lemma_1} and provided $\mu>0$, there exists $\lambda(t)>0$ such that the bilinear form 
    $a(q,\phi) = \int_{\Omega}(\mu \nabla q \cdot \nabla \phi-\mathbf{v}(t) \cdot \nabla \phi \, q) d \Omega$ satisfies the G\aa rding inequality
\begin{equation*}
    a(q,q)+\lambda(t) \norm{q}^2_{L^2(\Omega)} \geq \alpha_0(t) \norm{q}^2_{H^1(\Omega)}
\end{equation*}

for some $\alpha_0(t) \geq 0$. In particular, we can choose $\lambda(t)$ as
\begin{equation}
\label{lambda_eq}
    \lambda(t) = \frac{\norm{\mathbf{v}(t)}^2_{L^{\infty}(\Omega)^2}}{\mu}.
\end{equation}
We can also set:
\begin{equation}
\label{bar_alpha_0}
\bar{\alpha}_0 = \displaystyle \min_{0 \leq t \leq T }\left\{\frac{\mu}{2},\frac{\norm{\bb{v}(t)}^2_{L^{\infty}(\Omega)^2}}{2 \mu} \right\}.
\end{equation}
\qquad Proof: see the Appendix.
\QEDA
\end{lemma}

Now, we can prove a series of estimates on the norm of the state and its time derivative that ensure the well-posedness of the state equation, and will also be useful when showing the existence theorem for the OCP in Section \ref{OCP}.
\begin{theorem}[Well-posedness of the state problem ]
\label{theorem_1}
Assume that $\mathbf{u} \in \mathcal{U}_{ad} $, the initial density $q_0 \in L^2(\Omega) $ and $\mu >0$ is finite.
Then, there exists a unique weak solution $q \in L^2(0,T;H^1(\Omega))$ to the state problem (\ref{state_eq}) with $\dot{q} \in L^2(0,T;H^{1}(\Omega)^{*})$, such that the following energy estimates hold:
\begin{subequations}
\begin{align}
    &\norm{q}_{L^{\infty}(0,T;L^2(\Omega))}^2 \leq e^{\frac{16}{\mu} \norm{\mathbf{u}}^2_{\mathcal{U}}} \norm{q_0}_{L^2(\Omega)}^2, \label{subeq1} \\
    &\norm{q}_{L^{2}(0,T;L^2(\Omega))}^2 \leq T\,e^{\frac{16}{\mu} \norm{\mathbf{u}}_{\mathcal{U}}^2} \norm{q_0}_{L^2(\Omega)}^2, \label{subeq2}\\
    &\norm{q}^2_{L^2(0,T;H^1(\Omega))}  \leq \frac{1}{\bar{\alpha}_0}\Big(\frac{1}{2}+\frac{8}{\mu}\norm{\mathbf{u}}^2_{\mathcal{U}} e^{\frac{16}{\mu} \norm{\mathbf{u}}^2_{\mathcal{U}}}\Big)\norm{q_0}_{L^2(\Omega)}^2, \label{subeq3}\\
    &  \norm{\dot{q}}_{L^2(0,T;H^{1}(\Omega)^*)}^2  \notag \\ 
    &\leq \Big(\frac{\mu^2}{\bar{\alpha}_0}+ 16\Big(1+\frac{\mu}{\bar{\alpha}_0}\Big)\norm{\mathbf{u}}^2_{\mathcal{U}} \, e^{\frac{16}{\mu} \norm{\mathbf{u}}^2_{\mathcal{U}}}\Big)\norm{q_0}_{L^2(\Omega)}^2, \label{subeq4}
\end{align}
\end{subequations}
where $\bar{\alpha}_0$ is the minimum weak coercivity constant defined in Equation (\ref{bar_alpha_0}).

Proof: see the Appendix.
\QEDA

\end{theorem}

\section{The optimal control problem}
\label{OCP}

In this section, we state our OCP and derive a system of first-order necessary optimality conditions in the continuous framework, using the Lagrange multiplier approach, without any approximation on the form of the control functions and on the resulting state density function. Before doing that, we provide a mathematical analysis of the OCP by showing the existence of optimal controls and the differentiability of the control-to-state map.

The control problem can be framed as finding the boundary control actuation functions $\mathbf{u} \in \mathcal{U}_{ad}$ such that an initial density $q_0(\x)$ is optimally steered towards a target density $q_T(\x)$ while using as little actuation as possible. Obviously, these are conflicting objectives since zero actuation would end up in a uniform distribution due to diffusion. Formally, we can encode our objectives in a quadratic cost functional and write the OCP as:

\begin{equation}
\label{opt_problem_2}
\begin{aligned}
&  \tilde{J}(q,\mathbf{u}) \quad \longrightarrow \quad \min_{q,\bb{u}}   \\
&\textrm{s.t.}         \\
& \begin{cases}
 \quad &   \frac{\partial q}{\partial t} +  \nabla \cdot (-\mu\nabla q + \bb{v} q)   =0   \hspace{1cm}  in \,\, \Omega \times (0,T) \\
& (-\mu \nabla q + \bb{v} q)\cdot \nor = 0 \hspace{1.9cm} on \,\, \Gamma \times (0,T)
\phantom{space}  \\
&                         q(\x,0) = q_0(\x)  \hspace{2.8cm} in \,\, \Omega \,\, at \,\, t=0 \\
\end{cases} \\
& \textrm{and} \\
& \begin{cases}                           
                &\quad 0 \leq u_i(x_1,t) \leq u_{max}, \quad i=\{1,3\} \\
                 &  \quad 0 \leq u_i(x_2,t) \leq u_{max},\quad i=\{2,4\}. \\
                \end{cases}
\end{aligned}
\end{equation}
where 
\begin{equation*}
\begin{split}
\tilde{J}(q,\mathbf{u}) = \frac{1}{2} \int_{\Omega}  (q(\x,T)-q_T(\x))^2 \, d \Omega \\ + \frac{\alpha}{2} \sum_{i=1}^4  \int_{0}^T  \int_{ \Gamma_i} u_i(s,t)^2 \, d\Gamma\, dt.
\end{split}
\end{equation*}
\subsection{Analysis of the Optimal Control Problem}
In this subsection we prove the existence of an optimal control for Problem (\ref{opt_problem_2}) and the differentiability of the control-to-state map. For a fixed initial condition $q_0 \in L^{2}(\Omega)$, we define $q = \Xi(\mathbf{u})$ as the control-to-state-map, that is the state dynamics generated by the control function $\mathbf{u} \in \mathcal{U}_{ad}$.

\begin{theorem}[Differentiability of the control-to-state map]
\label{diff_c_t_s}
The control-to-state map $q = \Xi(\mathbf{u})$ is Fréchet differentiable and the directional derivative $z = \Xi'[\mathbf{u}]\mathbf{h}$ at $\mathbf{u} \in \mathcal{U}_{ad}$ in the direction $\mathbf{h} \in \mathcal{U}_{ad}$ is the solution of:
\begin{equation}
    \begin{cases}
\displaystyle \frac{\partial z}{\partial t}+\nabla \cdot(-\mu \nabla z+\bb{v}_{\bb{u}} z)=-\nabla \cdot(\bb{v}_{\bb{h}} q)  \hspace{0.3cm}  in \,\, \Omega \times (0,T) \\
(-\mu \nabla z+\bb{v}_{\bb{u}} z) \cdot \mathbf{n}=-\bb{v}_{\bb{h}} \cdot \mathbf{n} q \hspace{1.61cm} on \,\, \Gamma \times (0,T) \\
z(\mathbf{x}, 0)=0 \hspace{4.5cm} in \,\, \Omega \,\, at \,\, t=0
\end{cases}
\label{c_t_s_problem}
\tag{\ref{opt_problem_2} bis}
\end{equation}

where $\bb{v}_{\bb{u}}$ and $\bb{v}_{\bb{h}}$ are the velocity vector fields generated by $\bb{u}$ and $\bb{h}$ respectively, and $q = \Xi(\mathbf{u})$. 

 Proof: see the Appendix.
\QEDA
\end{theorem}

Finally, we show that at least an optimal control exists for the OCP (\ref{opt_problem_2}).

\begin{theorem}[Existence of an optimal control]
Let $q_0 \in L^2(\Omega)$. Consider the minimization problem of the reduced cost functional $J(\mathbf{u}) = \tilde{J}(\Xi(\mathbf{u}),\mathbf{u})$ over $\mathcal{U}_{ad}$, where $\tilde{J}$ is defined in (\ref{opt_problem_2}). Then, there exists a pair $(\bar{\mathbf{u}},\bar{q})$ such that $\bar{q}=\Xi(\bar{\mathbf{u}})$ and $\bar{\mathbf{u}}$ minimizes $J$ on $\mathcal{U}_{ad}$.

Proof: 

We already know that the state problem (\ref{opt_problem_2}) is well posed in the space $Y=H^{1}(0,T;H^{1}(\Omega),H^{1}(\Omega)^{*})$ where $H^{1}(0,T;H^{1}(\Omega),H^{1}(\Omega)^{*}) = \left\{ y \in L^2(0,T;H^{1}(\Omega)) : \dot{y} \in L^2(0,T;H^{1}(\Omega)^{*}) \right\}$, that is for every $\ut \in \mathcal{U}_{ad}$, problem (\ref{opt_problem_2}) has a unique solution $q=\Xi(\ut) \in Y$ and that, thanks to the estimates (\ref{subeq3}) and (\ref{subeq4}),
\begin{equation}
\label{est_Y}
    \norm{q}_{Y}^2 \leq C_0(\norm{\mbf{u}}_{\mathcal{U}}^2) \norm{q_0}_{L^2(\Omega)}^2
\end{equation}

\noindent where $C_0 = \Big(\frac{1}{\bar{\alpha}_0}\Big(\frac{1}{2}+\frac{8}{\mu}\norm{\mathbf{u}}^2_{\mathcal{U}} e^{\frac{16}{\mu} \norm{\mathbf{u}}^2_{\mathcal{U}}}\Big)+\Big(\frac{\mu^2}{\bar{\alpha}_0}+ 16\Big(1+\frac{\mu}{\bar{\alpha}_0}\Big)\norm{\mathbf{u}}^2_{\mathcal{U}} \, e^{\frac{16}{\mu} \norm{\mathbf{u}}^2_{\mathcal{U}}}\Big)\Big)$.

\noindent Moreover, the control-to-state map is $\mathcal{F}$-differentiable (see Theorem \ref{diff_c_t_s}). 

First of all, $\inf_{(q,\ut) \in Y \times \mathcal{U}} \tilde{J}(q,\ut)=I > -\infty$ and the set of feasible points is nonempty since $\tilde{J}(q,\ut)\geq 0 $.

For the sake of simplicity, we redefine the weak form of the state Equation (\ref{weak_form}) as: find $q$ such that $\forall \phi \in H^1(\Omega)$
\begin{equation}
    \langle \dot{q}(t),\phi\rangle_{H^{1}(\Omega)^{*},H^1(\Omega)} + d(q(t),\phi) + b(q(t),\phi;\mathbf{u}(t)) = 0,
\end{equation}
where $\langle \dot{q},\phi\rangle_{H^{1}(\Omega)^*,H^1(\Omega)} = \int_{\Omega} \dot{q} \phi \, d\Omega$.

We start by defining the operators associated to the bilinear forms $d$ and $b$ defined in Equation (\ref{bil_form}), then, using the estimates on the state equations we prove that they are bounded. The functional  $D: H^{1}(\Omega) \mapsto H^{1}(\Omega)^*$ associated to the bilinear form $d$ is:
\begin{equation*}
\langle D\,q(t),\phi \rangle_{H^{1}(\Omega)^*,H^1(\Omega)} = d(q(t),\phi) = \mu \int_{\Omega} \nabla q(t) \cdot \nabla \phi \, d\Omega
\end{equation*}
$\forall \phi \in H^1(\Omega)$.
From the definition of norm in the dual space $H^{1}(\Omega)^{*}$ we have:
\begin{equation*}
    \norm{ D\,q(t) }_{H^{1}(\Omega)^*} \leq \mu \norm{q(t)}_{H^{1}(\Omega)}
\end{equation*}
while squaring and integrating in time we get:
\begin{equation*}
\begin{split}
&\norm{ D\,q }_{L^2(0,T;H^{1}(\Omega)^*)}^2 \leq \mu^2 \norm{q}^2_{L^2(0,T;H^{1}(\Omega))} \\
&\leq   \frac{\mu^2}{\bar{\alpha}_0}\Big(\frac{1}{2}+\frac{8}{\mu}\norm{\mathbf{u}}^2_{\mathcal{U}} e^{\frac{16}{\mu} \norm{\mathbf{u}}^2_{\mathcal{U}}}\Big)\norm{q_0}_{L^2(\Omega)}^2 \leq C
\end{split}
\end{equation*}
for some constant $C\geq 0$, since $\norm{\mathbf{u}}^2_{\mathcal{U}}$ is bounded and $q_0 \in L^2(\Omega)$. Note that we have used the estimate obtained in Theorem \ref{theorem_1}. The functional $B: H^{1}(\Omega) \mapsto H^{1}(\Omega)^*$ associated to the bilinear form b is:
\begin{equation*}
\begin{split}
   & \langle B(\mathbf{u}(t),q(t)),\phi \rangle_{H^{1}(\Omega)^*,H^1(\Omega)} = b(q(t),\phi;\mathbf{u}(t)) \\
   &= \int_{\Omega} -q(t)\,\mathbf{v}(\mathbf{u}(t)) \cdot \nabla \phi \, d\Omega \qquad \forall \phi \in H^1(\Omega).
\end{split}
\end{equation*} 
Now we show that $B(\mathbf{u},q)$ is bounded. Using Cauchy-Schwarz inequality we obtain:
\begin{equation*}
\begin{split}
    &|\langle B(\mathbf{u}(t),q(t)),\phi \rangle_{H^{1}(\Omega)^{*},H^1(\Omega)}| = |b(q(t),\phi;\mathbf{u}(t))| \\
    &= \left|\int_{\Omega} q(t)\,\mathbf{v}(\mathbf{u}(t)) \cdot \nabla \phi \, d\Omega \right| \\
    &\leq \norm{\mathbf{v}(\mathbf{u}(t))}_{L^\infty(\Omega)^2} \norm{q(t)}_{L^2(\Omega)} \norm{\phi}_{H^1(\Omega)}
\end{split}
\end{equation*}
using the definition of norm in $H^{1}(\Omega)^{*}$ we have:
\begin{equation*}
    \norm{ B(\mathbf{u}(t),q(t))}_{H^{1}(\Omega)^{*}}\leq \norm{\mathbf{v}(\mathbf{u}(t))}_{L^\infty(\Omega)^2} \norm{q(t)}_{L^2(\Omega)},
\end{equation*}
while squaring and integrating over time we have:
\begin{equation*}
    \norm{ B(\mathbf{u},q)}_{L^2(0,T;H^{1}(\Omega)^{*})}^2 \leq \norm{q}_{L^{\infty}(0,T;L^2(\Omega))}^2 \norm{\mathbf{v}}_{L^2(0,T;L^\infty(\Omega)^2)}^2.
\end{equation*}
Finally, we use the estimates found in Lemma \ref{lemma_1} and Theorem \ref{theorem_1} to conclude that:
\begin{equation*}
    \norm{ B(\mathbf{u},q)}_{L^2(0,T;H^{1}(\Omega)^{*})}^2 \leq 8 e^{\frac{16}{\mu} \norm{\mathbf{u}}^2_{\mathcal{U}}} \norm{q_0}_{L^2(\Omega)}^2  \norm{\mathbf{u}}^2_{\mathcal{U}}  \leq C
\end{equation*}
for some constant $C\geq 0$. We can now recast the state equation in the dual space $H^{1}(\Omega)^{*}$ by defining:
\begin{equation}
\label{state_eq_abs}
G(q,\ut) = \begin{pmatrix}
   & \dot{q}(t) + Dq(t) + B(\mathbf{u}(t),q(t)) \\
   & q(0)-q_0
    \end{pmatrix}
\end{equation}
so that the state problem is $G(q,\ut)=0 \,\, a.e. \,\, t\in(0,T) $.

A minimizing sequence $\left\{(q_k,\ut_k) \right\}$ is bounded in $Y \times \mathcal{U}$.
Let $\{\mathbf{u}_n \}_{n\geq 1}$ be a minimizing sequence such that $\displaystyle \lim_{n\rightarrow \infty} J(\mathbf{u}_n) = I$ where $I = \displaystyle \inf_{\mathbf{u} \in \mathcal{U}_{ad}} J(\mathbf{u})$ and the associated sequence of states $\{ q_n \}_{n\geq 1}$ such that $q_n$ satisfies Equation (\ref{state_eq_abs}) for the control $\mathbf{u}_n$. Then, thanks to the definition of $\mathcal{U}_{ad}$, we have that $\norm{\mathbf{u}_n}_{\mathcal{U}}$ is bounded and from Equation (\ref{est_Y}) we deduce that $\norm{q_n}_{Y}$ is bounded as well.

The set of feasible points is weakly* sequentially closed in $Y \times \mathcal{U}$. Given the estimates on $q_n$ and $\mathbf{u}_n$, there exists a subsequence such that:

\begin{subequations}
\begin{align}
& \mathbf{u}_{n} \stackrel{\ast}{\rightharpoonup} \bar{\mathbf{u}} && \text { (weakly star) in } \mathcal{U} \notag \\
& q_{n} \stackrel{\ast}{\rightharpoonup} \bar{q} && \text { (weakly star) in } L^{\infty}(0, T ; L^2(\Omega)) \notag \\
& q_{n} \rightharpoonup \bar{q} && \text { (weakly) in } L^{2}\left(0, T ; H^1(\Omega)\right)  \label{est_contr_1}\\
& \dot{q}_{n} \rightharpoonup\psi && \text { (weakly) in } L^{2}\left(0, T ; H^{1}(\Omega)^{*}\right) \label{est_contr_2}\\
& D q_{n} \rightharpoonup \chi && \text { (weakly) in } L^{2}\left(0, T ; H^{1}(\Omega)^{*}\right) \notag \\
& B\left(\mathbf{u}_{n}, q_{n}\right) \rightharpoonup \Lambda && \text { (weakly) in } L^{2}\left(0, T ; H^{1}(\Omega)^{*}\right). \notag
\end{align}
\end{subequations}

\noindent Note that we have that $\psi = \dot{\bar{q}}$ and $\chi = D \bar{q}$. Moreover, from Equations (\ref{est_contr_1}) and (\ref{est_contr_2}) we have that $q_n \rightharpoonup \bar{q}$ in $Y$ and hence $q_n(0) \rightharpoonup \bar{q}(0) $ in $L^2(\Omega)$.

Finally we prove that $B(\bar{\mathbf{u}},\bar{q}) = \Lambda$. We can write:
\begin{equation}
\label{b_2_l}
\begin{aligned}
&\int_0^T \langle B(\bar{\bb{u}},\bar{q})-\Lambda,\phi \rangle \, dt   \\
&= \int_0^T\int_{\Omega} (-\bar{q}\,\mathbf{v}(\bar{\bb{u}}) + \lim_{n \rightarrow \infty} q_n\,\mathbf{v}(\bb{u}_n)) \cdot \nabla \phi \, d\Omega \, dt  \\
& = -\lim_{n \rightarrow \infty} \int_0^T\int_{\Omega} (\bar{q}\,\mathbf{v}(\bar{\bb{u}})   -q_n\,\mathbf{v}(\bb{u}_n)) \cdot \nabla \phi \, d\Omega d t \\
& = -\lim_{n \rightarrow \infty} \int_0^T\int_{\Omega} \bar{q} \, (\mathbf{v}(\bar{\bb{u}})   -\mathbf{v}(\bb{u}_n)) \cdot \nabla \phi \, d\Omega d t \\
&-\lim_{n \rightarrow \infty} \int_0^T\int_{\Omega} \mathbf{v}(\bb{u}_n)(\bar{q}-q_n) \cdot \nabla \phi \, d\Omega d t \qquad \forall \phi \in H^1(\Omega).
\end{aligned}
\end{equation}
Since $\bb{v}(\bb{u}) = \sum_{i=1}^4 \bb{b}_i u_i$ (see Appendix \ref{appendix_v}), we have $\bb{v}(\bar{\bb{u}})-\bb{v}(\bb{u}_n) = \sum_{i=1}^4 \bb{b}_i (\bar{u}_i - u_{i,n})$. 
Hence, for the first term in Equation (\ref{b_2_l}), we can write:
\begin{equation}
\label{limit_1}
\begin{aligned}
    & -\lim_{n \rightarrow \infty} \int_0^T\int_{\Omega} \bar{q} \, (\mathbf{v}(\bar{\bb{u}})   -\mathbf{v}(\bb{u}_n)) \cdot \nabla \phi \, d\Omega d t \\
    & =  \lim_{n \rightarrow \infty} \int_0^T\int_{\Omega} \bar{q} \, \sum_{i=1}^4 \bb{b}_i (u_{i,n}-\bar{u}_i) \cdot \nabla \phi \, d\Omega d t \\
    &= \lim_{n \rightarrow \infty} \sum_{i=1}^4  \int_0^T\int_{\Omega} \bar{q} \, \bb{b}_i\cdot \nabla \phi (u_{i,n}-\bar{u}_i)  \, d\Omega d t .
\end{aligned}
\end{equation}

We now prove that the limit in Equation \eqref{limit_1} is equal to zero. $\bb{b}_i$ is an analytic function in $\Omega$ thus $\bar{q}\,\bb{b}_i \cdot \nabla \phi  \in L^2(0,T;L^1(\Omega))$ and the domain $\Omega$ can be written as the Cartesian product $\Gamma_1 \times \Gamma_2 = (0,1) \times (0,1)$ therefore for the control function $u_1$ we have:
\begin{equation*}
\begin{aligned}
& \lim_{n \rightarrow \infty}   \int_0^T\int_{\Omega} \bar{q} \, \bb{b}_1\cdot \nabla \phi (u_{1,n}-\bar{u}_1)  \, d\Omega d t \\
& = \lim_{n \rightarrow \infty}   \int_0^T\int_{\Gamma_1}\int_{\Gamma_2} \bar{q} \, \bb{b}_1\cdot \nabla \phi  \,d\Gamma_2\,(u_{1,n}-\bar{u}_1)  \, d\Gamma_1 \, d t = 0
\end{aligned}
\end{equation*}
since, thanks to Fubini's Theorem, see e.g. \cite{brezis}, Chapter 4,  we have $\int_{\Gamma_2} \bar{q} \, \bb{b}_1\cdot \nabla \phi  \,d\Gamma_2 \in L^2(0,T;L^{1}(\Gamma_1))$ and $u_{1,n} \stackrel{\ast}{\rightharpoonup} \bar{u}_{1}$. The same resoning holds for the control functions $u_2,u_3,u_4$.

Note that, thanks to Aubin-Lions Lemma \cite{simon}, the embedding $H^1(0,T;H^1(\Omega),H^1(\Omega)^*) \hookrightarrow L^{2}(0, T ; L^{2}(\Omega))$ is compact. Thus, $\{q_n\}$ admits a subsequence strongly convergent to $\bar{q}$ in $L^2(0, T ; L^2(\Omega))$. Hence, regarding the second term in Equation (\ref{b_2_l}), we have:
\begin{equation*}
\begin{split}
     &\left|\int_0^T\int_{\Omega} \mathbf{v}(\bb{u}_n)(\bar{q}-q_n) \cdot \nabla \phi \, d\Omega d t \right| \\
     &\leq \norm{\mathbf{v}(\bb{u}_n)}_{L^{\infty}(Q)} \norm{\bar{q}-q_n}_{L^2(0,T;L^2(\Omega))} \norm{\nabla \phi}_{L^2(0,T;L^2(\Omega))},
\end{split}
\end{equation*}
so that:
\begin{equation*}
\lim_{n \rightarrow \infty} \left|\int_0^T\int_{\Omega} \mathbf{v}(\bb{u}_n)(\bar{q}-q_n) \cdot \nabla \phi \, d\Omega d t \right| = 0,
\end{equation*}

 \noindent since $q_n \rightarrow \bar{q}$ strongly  and $\norm{\mathbf{v}(\bb{u}_n)}_{L^{\infty}(Q)}$ is uniformly bounded.

 Since $\Gamma_i$ is bounded for $i=1 \ldots 4$, the weak star convergence of $u_{n,i}$ in $L^2(0,T;L^{\infty}(\Gamma_i))$  to some $\bar{u}_i \in L^2(0,T;L^{\infty}(\Gamma_i))$ implies weak convergence of $\bb{u}_n$ to $\bar{\bb{u}}$ in any $L^2(0,T;L^p(\Gamma))$, $1 \leq p < \infty$, and in particular in $L^2(0,T;L^2(\Gamma))$. Then, exploiting the fact that $q_n$ weakly converges to $\bar{q}$ in $L^2(0,T;H^1(\Omega))$ and that $\tilde{J}(q,{\bf u})$ is convex and continuous in $L^2(0,T;H^1(\Omega)) \times L^2(0,T;L^2(\Gamma))$, we conclude that:
\begin{equation*}
    J(\bar{\mathbf{u}}) \leq \displaystyle \lim_{n \rightarrow \infty} \inf J(\mathbf{u}_n) = I;
\end{equation*}
thus, the pair $(\bar{\bb{u}},\bar{q})$ is an optimal pair for the considered optimal control problem.
\QEDA
\end{theorem}

\subsection{Optimality conditions}
We now derive a set of first-order optimality conditions using the Lagrangian method \cite{fredi}. Using this idea, we obtain an explicit expression for the gradient of the cost functional in the continuous setting. The Lagrangian functional $\mathcal{L}: \mathcal{V} \times \mathcal{U} \times \mathcal{W}^{*} \mapsto \R$ is defined as

\begin{equation*}
    \mathcal{L}(q,\bb{u},p) = \tilde{J}(\bb{u},q) + \langle p , G(q,\bb{u}) \rangle_{\mathcal{W}^*,\mathcal{W}}
\end{equation*}
where $\mathcal{V}=H^{1}(0,T;H^{1}(\Omega),H^{1}(\Omega)^*)$ and $\mathcal{W}=L^2(0,T;H^1(\Omega)^{*})$ so that state, control and adjoint variables are considered independently. Therefore, the set of first-order necessary optimality conditions consists of imposing that the Gateaux derivative of the Lagrangian with respect to the triple $(q,\bb{u},p)$ along an arbitrary variation $(\psi,\bb{h},\phi)$ is equal to zero. 
In our case, the Lagrangian takes the explicit form
\begin{equation}
\label{lag_def}
    \mathcal{L} = \tilde{J}(\bb{u},q) - \int_{\Omega}\int_0^T \frac{\partial q}{\partial t}p +  \nabla \cdot (-\mu\nabla q + \bb{v} q)p \,\, d \Omega d t ,
\end{equation}
where $p \colon \Omega \times [0,T] \mapsto \R$ is the adjoint variable relative to the dynamic constraint. In order to derive the adjoint system fulfilled by $p$, it is useful to rearrange some terms of the Lagrangian so that we can express the portion of the Lagrangian relative to the dynamic constraint as:
\begin{equation*}
\begin{aligned}
&\int_{\Omega}\int_0^T \frac{\partial q}{\partial t}p +  \nabla \cdot (-\mu\nabla q + \bb{v} q)p \,\, d \Omega d t \\
& =  \int_{\Omega} q(\bb{x},T)p(\bb{x},T) \, d\Omega  -\int_{\Omega}q(\bb{x},0)p(\bb{x},0) \, d\Omega \\
&- \int_{\Omega}\int_0^T q\frac{\partial p}{\partial t} \, dt d\Omega + \int_{\Gamma}\int_0^T (\mu \nabla  p q)\cdot \nor dt d\Gamma \\
&+\int_{\Omega}\int_0^T-\mu \Delta p  q \, dt d\Omega -\int_{\Omega}\int_0^T\bb{v} \cdot \nabla p \,\, dt d\Omega.
\end{aligned}
\end{equation*}

The adjoint dynamics is obtained by imposing $\mathcal{L}^{'}_{q}[\psi] = 0$, hence 
\begin{equation*}
\begin{aligned}
    &\mathcal{L}^{'}_{q}[\psi] =  \int_{\Omega}  (q(\x,T)-q_T(\x)) \,\psi(\x,T) \, d \Omega  \\
    &- \int_{\Omega}\int_0^T \left( -\mu \Delta p - \frac{\partial p}{\partial t} - \bb{v} \cdot \nabla p \right)\,\psi \, d \Omega \, dt  \\
    & -\int_{\Gamma}\int_0^T \mu \nabla  p \cdot \nor \, \psi \, dt d\Gamma \\
    & -\int_{\Omega} \psi(\x,T)p(\bb{x},T) \, d\Omega = 0
    \end{aligned}
\end{equation*}
where $\psi(\x,0)=0$ since the initial condition on $q(\x,t)$ is fixed. The adjoint dynamics thus reads:

\begin{equation}
\begin{aligned}
\label{adj_eq}
    &- \frac{\partial p}{\partial t}-\mu \Delta p  - \bb{v} \cdot \nabla p = 0  \hspace{1cm}  in \,\, \Omega \times (0,T)\\
    & \nabla  p \cdot \nor = 0                                                 \hspace{3.43cm} on \,\, \Gamma \times (0,T) \\
    & p(\x,T) = q(\x,T)-q_T(\x)                                                \hspace{1.1cm} in \,\, \Omega \,\, at \,\, t=T.\\ 
    \end{aligned}
\end{equation}
 Note that the adjoint problem is backward in time, and a final time condition is imposed. We also highlight that the velocity field $\bb{v}$ appears with an opposite sign compared to the state equation. Furthermore, even if the velocity field varies in space (that is $\nabla \cdot \bb{v} \neq 0$), this dependence does not directly affect the adjoint equation. Finally, the boundary conditions of the adjoint problem are of homogeneous Neumann type while in the state equation we had Neumann no-flux boundary conditions. Up to now, we have not taken into account the explicit dependence of the velocity field $\bb{v}$ from the set of control actions (i.e. the dependence $\bb{v}=\bb{v}(\bb{u})$) that is formalized in Equation (\ref{act}). Interestingly, this dependence does not affect the derivation of the adjoint equation.
 The reduced gradient is obtained by imposing $\mathcal{L}^{'}_{\bb{u}}[\bb{h}] = \bb{0}$. The gradient relative to each control function is derived by taking a variation along that direction only. For the sake of clarity, we rewrite the terms of the Lagrangian that depends on the control functions as:
\begin{equation*}
\begin{aligned}
    &\mathcal{L} = \ldots + \frac{\alpha}{2} \int_{T} \sum_{i=1}^4 \int_{\Gamma_i} u_i(x_i,t)^2 \, d\Gamma\, dt\\ &+\int_{\Omega}\int_0^T \bb{v} \cdot \nabla p \, q d\Omega dt  \\
    &=\ldots + \frac{\alpha}{2} \int_{T} \sum_{i=1}^4 \int_{\Gamma_i} u_i(x_i,t)^2 \, d\Gamma\, dt\\
    & +\int_{\Omega}\int_0^T q\,v_{1} \frac{\partial p}{\partial x_1}+q\,v_{2}\frac{\partial p}{\partial x_2} \, d\Omega dt
    \end{aligned}
\end{equation*}    
where, substituting the form of $v_{1}$ and $v_{2}$ in Equation (\ref{act}), we have:
\begin{equation*}
    \begin{aligned}
    &\int_{\Omega}\int_0^T q\,v_{1} \frac{\partial p}{\partial x_1}+q\,v_{2}\frac{\partial p}{\partial x_2} \, d\Omega dt  \\
    & = \int_{\Omega}\int_0^T q\,\Big( u_4(x_2,t) e^{-cx_1} - u_2(x_2,t) e^{-c(1-x_1)} \Big) \frac{\partial p}{\partial x_1}d\Omega dt\,\\
    & + \int_{\Omega}\int_0^T q \,\Big( u_1(x_1,t) e^{-cx_2} - u_3(x_1,t) e^{-c(1-x_2)} \Big) \frac{\partial p}{\partial x_2} \, d\Omega dt. \\
    \end{aligned}
\end{equation*}

\noindent The Gateaux derivative with respect to the control function $u_1$ along the direction $h_1$ is:
\begin{equation}
    \label{lag_cont}
    \begin{aligned}
    &\mathcal{L}^{'}_{u_1}[h_1] = \alpha \int_{0}^T \int_{\Gamma_1} u_1(x_1,t) h_1(x_1,t) \, dt \, d\Gamma \\
    &+ \int_{\Omega}\int_0^T e^{-c x_2} \, q\, \frac{\partial p}{\partial x_2} h_1(x_1,t) d\Omega dt.
    \end{aligned}
\end{equation}
In our case, $\Omega =(0,1)^2$ and the function $u_1(x_1,t): [0,1]\times[0,T] \rightarrow \R$ is a function of the $x_1$ variable only, thus, the integral in (\ref{lag_cont}) can be split and simplified as:
\begin{equation*}
\begin{aligned}
&\mathcal{L}^{'}_{u_1}[h_1] =  \alpha \int_{0}^T \int_{0}^1 u_1 \, h_1 \, dx_1 \, dt \,   \\
    &+\int_{0}^1\int_{0}^1\int_0^T e^{-cx_2} \, q\, \frac{\partial p}{\partial x_2} h_1 \,\,dx_1 \, dx_2 \, dt  \\
    & =\int_{0}^T\int_{0}^1 \left[\alpha u_1 +  \int_{0}^1 \, q\, e^{-c x_2} \frac{\partial p}{\partial x_2} \, dx_2\right] \,  h_1 \, dx_1 \, dt \\
    &=\int_{0}^T\int_{0}^1 \nabla J_1(x_1,t) h_1(x_1,t) \,\, dx_1 \, dt \, ,
\end{aligned}
\end{equation*}
\noindent for any variation of the first control function $h_1(x_1,t)$. 
We have identified the reduced gradient of the cost functional $J$ with respect to the control function $u_1$ as:
\begin{equation*}
    \nabla J_1(x_1,t) = \alpha u_1(x_1,t) +  \int_{0}^1 \, q\, e^{-cx_2} \frac{\partial p(x_1,x_2,t)}{\partial x_2} \, dx_2.
\end{equation*}
The gradients with respect to $u_2,u_3$ and $u_4$ are obtained in a similar way.
In the continuous setting, the set of first-order necessary conditions that the optimal triple $(q^{\star},\ut^{\star},p^{\star})$ must satisfy is given by:
\begin{itemize}

\item the optimal state dynamics: 
\begin{equation*}
\begin{aligned}
    &\frac{\partial q^{\star}}{\partial t} +  \nabla \cdot (-\mu\nabla q^{\star} + \bb{v}_{\bb{u}^{\star}} q^{\star})   =0  \hspace{0.5cm}  in \,\, \Omega \times (0,T)\\
    &(-\mu \nabla q^{\star} + \bb{v}_{\bb{u}^{\star}} q^{\star})\cdot \nor = 0 \hspace{1.63cm} on \,\, \Gamma \times (0,T) \\
    & q^{\star}(\x,0) = q_0(\x) \hspace{3.05cm} in \,\, \Omega \,\, at \,\, t=0\\ 
    \end{aligned}
\end{equation*}

\item the adjoint dynamics: 

\begin{equation}
\begin{aligned}
\label{adj_eq}
&-\frac{\partial p^{\star}}{\partial t}-\mu \Delta p^{\star}  - \bb{v}_{\bb{u}^{\star}} \cdot \nabla p^{\star} = 0 \hspace{0.8cm}  in \,\, \Omega \times (0,T)\\
&\nabla  p^{\star} \cdot \nor = 0  \hspace{3.85cm} on \,\, \Gamma \times (0,T) \\
&p^{\star}(\x,T) = q^{\star}(\x,T)-q_T(\x)  \hspace{1.4cm} in \,\, \Omega \,\, at \,\, t=T\\ \end{aligned}
\end{equation}

\end{itemize}

and the variational inequalities:

\begin{equation*}
\int_0^1 \int_0^T \nabla J_i  (u_i-u_i^{\star}) \,\, dx_k \, dt \geq 0 \quad \forall i=1,\ldots,4 \quad \forall \bb{u} \in \mathcal{U}_{ad}
\end{equation*}
where $k=1$ for $i=\{1,3\}$ and $k=2$ for $i=\{2,4\}$. 
Note that the vector field $\bb{v}_{\bb{u}}$ depends on the control vector as shown in Equation (\ref{act}); see the Appendix \ref{appendix_v} for further details. The functional relationship $\bb{v}=\bb{v}(\bb{u})$ affects the expression of the gradient of the reduced cost functional while does not influence directly the form of the adjoint problem.

\begin{remark}[Controllability issues]
The controllability properties of the Kolmogorov forward equation with no-flux boundary conditions have been studied in \cite{karthik_bil} obtaining weaker conditions on the target and initial densities. \cite{Duprez2019} considers the controllability problem around a nominal trajectory where the control vector field is defined in a compact region of the workspace, as stated in the paper their derivation strongly relies on the Dirichlet boundary conditions and results concerning a no-flux boundary are yet to be established. Very recently, \cite{noflux_barbu} established asymptotic controllability with a scalar control in the diffusive term and no-flux boundary conditions. In these works the control vector field is the velocity field itself, either in the whole domain or in a compact subset of it. In our case, the control actions are defined along the boundary of the domain and are mapped to the velocity field though a function that depends on the actuation mechanism. The rigorous controllability analysis for our case is left to future work but we show through numerical test cases that relatively complex target densities can be approximately reached. As can be intuitively guessed, closed shapes cannot be reached with sufficient accuracy.
\end{remark}

\section{Numerical approximation}
\label{Num_approx}
 In this Section, the control functions are approximated with one-dimensional Radial Basis Functions (RBF) while the PDEs arising from state and adjoint dynamics are discretized in space with the Finite Element Method (FEM) thus pursuing an ``Optimize-then-Discretize" (OtD) approach.
 Afterward, we reformulate the control problem in the so-called ``Discretize then Optimize" (DtO) framework and we draw some similarities between the resulting control problems obtained following either strategies.

\subsection{Optimize-then-Discretize}

\noindent The space of continuous functions in which the set of control functions is taken is approximated in space with a finite number of basis functions as:
\begin{equation}
    \label{control_approx}
    u_1(x_1,t) = \sum_{i=1}^{N_c} \psi_i(x_1) u_{1,i}(t);
\end{equation}
under this choice applying Theorem \ref{theorem_1} is straightforward if we take an appropriate set of basis functions $\psi_i(x_1)$.
The form of the adjoint Equation (\ref{adj_eq}) is unchanged, we just need to specialize the form of the vector field $\bb{v}$ determined by the approximation of control functions $u_k$ in its expression as:
\begin{equation*}
    \begin{aligned}
    &- \frac{\partial p}{\partial t}-\mu \Delta p  - \bb{v} \cdot \nabla p  = - \frac{\partial p}{\partial t}-\mu \Delta p   \\
    &-\sum_{i=1}^{N_c} \psi_i(x_2) \Big(u_{4,i}(t) e^{-cx_1} - u_{2,i}(t) e^{-c(1-x_1)}\Big) \frac{\partial p}{\partial x_1}  \\
    &-\sum_{i=1}^{N_c} \psi_i(x_1)\Big( u_{1,i}(t) e^{-cx_2} - u_{3,i}(t)  e^{-c(1-x_2)}\Big) \frac{\partial p}{\partial x_2} = 0,
\end{aligned}
\end{equation*}

whereas to recover the form of the gradient it is easier to modify the term in the cost functional that weights the control actions as:
\begin{equation*}
    J_c = \frac{\alpha}{2}\int_0^T \sum_{k=1}^4 \sum_{i=1}^{N_c} u_{k,i}^2(t) \, dt.
\end{equation*}
This change makes sense since we have no control on the $\psi_i$ basis functions that had been chosen beforehand. The reduced gradient is recovered by taking a variation along each $u_{1,i}$. We will derive the optimality conditions for the generic set of basis functions $\psi_i$ and then adapt it to a specific choice later. The reduced gradient for the control coefficient related to the $i$-th basis function of actuator $1$ is:

\begin{equation*}
    \nabla J_{1,i}(t) = \alpha u_{1,i}(t) + \int_{\Omega} \psi_i(x_1) e^{-cx_2} \, q \, \frac{\partial p}{\partial x_2} \, d \Omega,
\end{equation*}
while the other components are obtained in the same way.
Note that since we have prescribed a fixed spatial shape for the control functions, the reduced gradients are functions of time only and measure the sensitivity to variations in the control coefficients $u_{k,i}$ for the $i$-th coefficients of the $k$-th control function. Since there are $4$ control functions (i.e. one for each side) and $N_c$ basis functions, we have $4 \times N_c$ reduced gradients as functions of time.

We proceed in the numerical approximation of the OCP resorting to a Finite Element Method (FEM). The state and adjoint variables are expressed as linear combinations of a set of basis functions, using the same basis for both state and adjoint problems. This choice allows us to draw an interesting comparison between the OtD and DtO approaches. The state and adjoint functions are then approximated as:
\begin{equation}
\label{state_approx}
    q(\x,t) = \sum_{i=1}^{N_q} \phi_i(\x)q_i(t),
\end{equation}
\begin{equation}
\label{adjoint_approx}
    p(\x,t) = \sum_{i=1}^{N_p} \phi_i(\x)p_i(t) ,
\end{equation}
\noindent where we also choose $N_p=N_q$. We select piecewise linear, globally continuous ansatz functions $\phi_i$ ($\mathbb{P}_1$ finite elements). The derivation of the resulting system of ordinary differential equations (ODEs) is standard practice in the treatment of linear parabolic PDEs with the FEM. However, some care is needed for the case at hand because the control functions appear in the Neumann boundary conditions of the state equation. Thus, we carry out the complete derivation of the ODE system starting from the weak formulation of the state equation (\ref{state_weak_2}). 

First of all, we substitute in Equation (\ref{state_weak_2}) the state approximation (\ref{state_approx}) and the form of $\bb{v}$ in Equation (\ref{act}) with the spatial control approximation in Equation (\ref{control_approx}). We rename $w$ the test function and we assume it belongs to the same finite dimensional space, that is $w=\sum_{i=1}^{N_q} \phi_i\,w_i$. The state equation (\ref{state_weak_2}) must be satisfied for each basis function $\phi_i$.
Carrying out the substitutions we obtain the usual mass matrix $M$, whose elements are defined as:

\begin{equation*}
    (M)_{ij} = \int_{\Omega} \phi_j(\x) \, \phi_i(\x)\, d \Omega,
\end{equation*}

and the diffusion matrix $A$ whose elements are:

\begin{equation*}
    (A)_{ij} = \int_{\Omega} \mu \, \nabla \phi_j(\x) \cdot \nabla\phi_i(\x) \, d\Omega.
\end{equation*}

The pure transport term is somewhat more involved and need to be split as:
\begin{equation*}
    \int_{\Omega}  \bb{v}\cdot \nabla q \, \phi_i \, d\Omega =\int_{\Omega}  v_{1} \frac{\partial q}{\partial x_{1}} \, \phi_i \, d\Omega +\int_{\Omega}  v_{2} \frac{\partial q}{\partial x_2 } \, \phi_i \, d\Omega,
\end{equation*}
for the $v_{1}$ term we have:
\begin{equation*}
\begin{aligned}
    &\int_{\Omega}  v_{1} \frac{\partial q}{\partial x_1} \, \phi_i \, d\Omega \\
    &= \int_{\Omega}  \Big(u_4(x_2,t) e^{-cx_1} - u_2(x_2,t) e^{-c(1-x_1)} \Big) \frac{\partial q}{\partial x_1} \, \phi_i \, d\Omega \\
    &= \int_{\Omega} u_4(x_2,t) e^{-c x_1} \frac{\partial q}{\partial x_1}  \phi_i \, d\Omega  \\
    &- \int_{\Omega} u_2(x_2,t) e^{-c(1-x_1)} \frac{\partial q}{\partial x_1}  \phi_i \, d\Omega  \\
    & =\Big(\int_{\Omega} \sum_{k=1}^{N_c} \sum_{j=1}^{N_q} \psi_k \, e^{-c x_1} \, \frac{\partial \phi_j}{\partial x_1} \, \phi_i \, d\Omega \Big) \,  u_{4,k}(t) q_j(t)  \\
    & - \Big(\int_{\Omega} \sum_{k=1}^{N_c} \sum_{j=1}^{N_q} \psi_k \, e^{-c(1-x_1)} \, \frac{\partial \phi_j}{\partial x_1} \, \phi_i \, d\Omega \Big) \,  u_{2,k}(t) q_j(t).
\end{aligned}
\end{equation*}
We define the matrices $B_{4,k}$ and $B_{2,k}$ as:
\begin{equation*}
(B_{4,k})_{ij} = \int_{\Omega} \psi_k(x_2) \, e^{-cx_1} \, \frac{\partial \phi_j(\x)}{\partial x_1} \, \phi_i(\x) \, d\Omega,
\end{equation*}
\begin{equation*}
(B_{2,k})_{ij} = \int_{\Omega} \psi_k(x_2) \, e^{-c(1-x_1)} \, \frac{\partial \phi_j(\x)}{\partial x_1} \, \phi_i(\x) \, d\Omega.
\end{equation*}
Carrying out the same substitutions, the $v_{2}$ component determines the matrices $B_{1,k}$ and $B_{3,k}$.

Now we examine the reaction term generated by the divergence of the velocity vector field:
\begin{equation*}
    \int_{\Omega}  (\nabla \cdot \bb{v}) \, q \, \phi_i \,\, d\Omega = \int_{\Omega}  \frac{\partial v_1}{\partial x_1} \, q \, \phi_i \,\, d\Omega + \int_{\Omega}  \frac{\partial v_2}{\partial x_2} \, q \, \phi_i \,\, d\Omega;
\end{equation*}
regarding the $v_{1}$ term we note that:
\begin{equation*}
    \frac{\partial v_{1}}{\partial x_1} = \frac{\partial }{\partial x_1} \Big(u_4(x_2,t) e^{-cx_1} - u_2(x_2,t) e^{-c(1-x_1)} \Big) = - c v_{1},
\end{equation*}
and thus:
\begin{equation*}
\begin{aligned}
&\int_{\Omega}  \frac{\partial v_{1}}{\partial x_1} \, q \, \phi_i \,\, d\Omega  \\
&= \int_{\Omega}  - c \, v_{1} \, q \, \phi_i \,\, d\Omega   \\
&= \int_{\Omega}  - c \Big(u_4(x_2,t) e^{-cx_1} - u_2(x_2,t) e^{-c(1-x_1)}\Big) \, q \, \phi_i \,\, d\Omega   \\ 
    & = \Big(\int_{\Omega} \sum_{k=1}^{N_c} \sum_{j=1}^{N_q} -c \psi_k \,  e^{-cx_1} \, \phi_j \, \phi_i \, d\Omega \Big) \,  u_{4,k}(t) q_j(t)  \\
& - \Big(\int_{\Omega} \sum_{k=1}^{N_c} \sum_{j=1}^{N_q} - c \psi_k \,  e^{-c(1-x_1)} \, \phi_j \, \phi_i \, d\Omega \Big) \,  u_{2,k}(t) q_j(t).
\end{aligned}
\end{equation*}
We define the matrices $C_{4,k}$ and $C_{2,k}$ as:
\begin{equation*}
    (C_{4,k})_{ij} = \int_{\Omega}  -c \psi_k(x_2) \,  e^{-cx_1} \, \phi_j(\x) \, \phi_i(\x) \, d\Omega 
\end{equation*}
\begin{equation*}
    (C_{2,k})_{ij} = \int_{\Omega}  -c \psi_k(x_2) \,  e^{-c(1-x_1)} \, \phi_j(\x) \, \phi_i(\x) \, d\Omega.
\end{equation*}
With the same idea we carry out the derivation of the matrices $C_{1,k}$ and $C_{3,k}$ regarding the $v_{2}$ terms.

Finally, we treat the term generated by the boundary conditions:
\begin{equation*}
\begin{aligned}
&\int_{\Gamma} -\bb{v} q \cdot \bb{n} \, \phi_i \, d\Gamma   \\
&=\int_{\Gamma} -v_{1} \,n_{1} \, q  \, \phi_i \, d\Gamma + \int_{\Gamma} -v_{2}\, n_{2}\, q  \, \phi_i \, d\Gamma ;
\end{aligned}
\end{equation*}
regarding the $v_{1}$ term we have:
\begin{equation*}
\begin{aligned}
&\int_{\Gamma} -v_{1} n_{1} \, q \, w \, d\Gamma  \\
&=\int_{\Gamma} - \Big(u_4(x_2,t) e^{-cx_1} - u_2(x_2,t) e^{-c(1-x_1)} \Big) \, n_1 \, q \, w \, d\Gamma  \\
& = \Big(\int_{\Gamma}  \sum_{k=1}^{N_c}\sum_{j=1}^{N_q} -\psi_k(x_2) e^{-cx_1} \phi_j \, \phi_i n_{1} \, d \Gamma\Big) u_{4,k}(t) q_{j}(t)  \\
&-\Big(\int_{\Gamma}  \sum_{k=1}^{N_c}\sum_{j=1}^{N_q} -\psi_k(x_2) e^{-c(1-x_1)} \phi_j \, \phi_i n_{1}\, d \Gamma\Big) u_{2,k}(t) q_{j}(t),
\end{aligned}
\end{equation*}
and we define matrices $L_{4,k}$ and $L_{2,k}$ as:
\begin{equation}
\label{L_4}
    (L_{4,k})_{ij} = \int_{\Gamma} -\psi_k(x_2) e^{-c{x_1}} \phi_j(\x) \, \phi_i(\x) n_{1}(\x) \, d \Gamma
\end{equation}
\begin{equation*}
        (L_{2,k})_{ij} = \int_{\Gamma} -\psi_k(x_2) e^{-c(1-x_1)} \phi_j(\x) \, \phi_i(\x) n_{1}(\x) \, d \Gamma.
\end{equation*}
Similarly, for the $v_{2}$ term we define the matrices $L_{1,k}$ and $L_{3,k}$.

\noindent Grouping the unknown time-dependent coefficients of the state in a vector $\bb{q} = [q_1,\ldots,q_{N_q}]^{\top}$ we obtain the following system of ODEs arising from the spatial discretization of state and control functions:
\begin{equation*}
    \begin{cases}
    &M \dot{\mbf{q}} + A \mbf{q} + \Gamma_q(\ut) \mbf{q} = \mbf{0}  \quad , \quad t \in (0,T)\\
    & \mbf{q}(0)=\bar{\mbf{q}}_0
    \end{cases} 
\end{equation*}
where $\Gamma_q(\ut) \in \R^{N_q \times N_q}$ is defined as 
\begin{equation}
\label{gamma_q}
\begin{aligned}
    &\Gamma_q(\ut) = \sum_{k=1}^{N_c} \Big\{u_{4,k}\Big(B_{4,k}+C_{4,k}+L_{4,k}\Big) \\ &+u_{1,k}\Big(B_{1,k}+C_{1,k}+L_{1,k}\Big)  \\
    & - u_{3,k}\Big(B_{3,k}+C_{3,k}+L_{3,k}\Big)  \\
    & -u_{2,k}\Big(B_{2,k}+C_{2,k}+L_{2,k}\Big)\Big\}
\end{aligned}
\end{equation}
that is, a bilinear control system with state $\mbf{q}(t) \in \R^{N_q}$ and control inputs $u_{c,k}(t)$, $c=1,\ldots,4$ and $k=1,\ldots,N_c$. The discretization (in space) of the problem data is denoted by an overline (e.g. $\bar{\mbf{q}}_0$ corresponds to the given initial condition).

\begin{lemma}[Properties of the B,C,L matrices]
\label{BCL}
The matrices $B_{c,k}$,\,\,$C_{c,k}$,\,\,$L_{c,k}$  satisfy:
\begin{equation}
B_{c,k}+B^{\top}_{c,k}+C_{c,k}+L_{c,k} = 0
\end{equation}
for any choice of FEM ansatz functions $\phi_i$, control basis functions $\psi_k$ $k=1,\ldots,N_c$ and actuator side $c=1,\ldots,4$.

Proof:
We will just prove it for $L_{4,k}$. The results for all the other matrices are identical. Recalling the definition of $(L_{4,k})_{ij}$ in Equation (\ref{L_4}) and defining $\tilde{\psi}_k(x_1,x_2) :=  \psi_k(x_2) e^{-cx_1}$, using integration by parts we get:

\begin{equation*}
\begin{aligned}
        &(L_{4,k})_{ij} = \int_{\Gamma} -\psi_k(x_2) e^{-c x_1} \phi_j(\x) \, \phi_i(\x) n_{1}(\x) \, d \Gamma   \\
        &= -\int_{\Gamma} \tilde{\psi}_k(\x) \phi_j(\x) \, \phi_i(\x) n_{1}(\x) \, d \Gamma  \\
        & -\int_{ \Omega} \frac{\partial}{\partial x_1}\Big(\tilde{\psi}_k(\x) \phi_j(\x) \, \phi_i(\x)\Big) \, d \Omega  \\
        & -\int_{ \Omega} \frac{\partial \phi_i}{\partial x_1 }  \tilde{\psi}_k \phi_j  \, d\Omega - \int_{ \Omega} \frac{\partial \phi_j}{\partial x_1 }  \tilde{\psi}_k \phi_i  \, d\Omega  \\
        & -\int_{ \Omega} \frac{\partial \tilde{\psi}_k}{\partial x_1 }   \phi_i \phi_j  \, d\Omega = -(B_{4,k})_{ji}-(B_{4,k})_{ij}-(C_{4,k})_{ij}
  \end{aligned}
\end{equation*}
and thus: 
\begin{equation*}
    (L_{4,k})_{ij} + (B_{4,k})_{ji}+(B_{4,k})_{ij}+(C_{4,k})_{ij} = 0.
\end{equation*}
\QEDA
\end{lemma}

With the same reasoning and definitions for the state and input matrices the FEM discretization of the adjoint equation reads:

\begin{equation}
\label{otd_adjoint}
    \begin{cases}
    & -M \dot{\mbf{p}} + A \mbf{p} + \Gamma_p(\ut) \mbf{p} = \mbf{0}\quad , \quad t \in (0,T)  \\
    &\mbf{p}(T) = \mbf{q}(T)-\bar{\mbf{q}}_T
    \end{cases}
\end{equation}
where the matrix $\Gamma_p \in \R^{N_q \times N_q}$ is defined as: 
\begin{equation}
\begin{aligned}
\label{gamma_p}
    &\Gamma_p(\ut) = -\sum_{k=1}^{N_c} \Big\{ u_{4,k}B_{4,k}+u_{1,k}B_{1,k}  \\
    & - u_{3,k}B_{3,k}-u_{2,k}B_{2,k} \Big\}.
\end{aligned}
\end{equation}
Note that the matrices $C_{a,k}$ and $L_{a,k}$ do not appear in the adjoint equations while the matrices $B_{a,k}$ related to pure transport terms appear with opposite sign. Note also that, due to Lemma \ref{BCL}, we have $-B_{a,k} = B^{\top}_{a,k}+C_{a,k}+L_{a,k}$, $C_{a,k}$ and $L_{a,k}$ being symmetric. It is possible to show that: 
\begin{equation}
\label{gamma_qp}
    \Gamma_q(\ut)^{\top} = \Gamma_p(\ut).
\end{equation}
From its definition in Equation (\ref{gamma_q}) we get:
\begin{equation*}
\begin{aligned}
&\Gamma_q(\ut)^{\top}  = \sum_{k=1}^{N_c} \Big\{u_{4,k}\Big(B^{\top}_{4,k}+C^{\top}_{4,k}+L^{\top}_{4,k}\Big) \\ &+u_{1,k}\Big(B^{\top}_{1,k}+C^{\top}_{1,k}+L^{\top}_{1,k}\Big)  \\
    & - u_{3,k}\Big(B^{\top}_{3,k}+C^{\top}_{3,k}+L^{\top}_{3,k}\Big)  \\
    & -u_{2,k}\Big(B^{\top}_{2,k}+C^{\top}_{2,k}+L^{\top}_{2,k}\Big)\Big\}  \\
    &= -\sum_{k=1}^{N_c} \Big\{ u_{4,k}B_{4,k}+u_{1,k}B_{1,k}  \\
    & - u_{3,k}B_{3,k}-u_{2,k}B_{2,k} \Big\} = \Gamma_p(\ut),
    \end{aligned}
    \end{equation*}
where we have used the symmetry of the matrices $C$ and $L$ and Lemma \ref{BCL}. Equation (\ref{gamma_qp}) will be useful to investigate the interconnection between OtD and DtO approaches.

\noindent The gradient can be expressed as:
\begin{equation}
\label{OtD_grad}
    \begin{aligned}
    & \nabla J_{1,k}(t) =\alpha \, u_{1,k}(t) + \int_{\Omega} \psi_k(x_1) e^{-c x_2} \, q \, \frac{\partial p}{\partial x_2} \, d \Omega  \\
    &=\alpha \, u_{1,k}(t) + \int_{\Omega} \psi_k(x_1) e^{-cx_2} \, \sum_{i=1}^{N_q}\sum_{j=1}^{N_p} \phi_i \, \frac{\partial \phi_j}{\partial x_2} q_i p_j \, d \Omega  \\
    &= \alpha \, u_{1,k}(t) + \mbf{q}(t)^{\top} \, B_{1,k} \, \mbf{p}(t),
    \end{aligned}
\end{equation}
where we have replaced the definition of $B_{1,k}$. $\nabla J_{2,k}(t)$,$\nabla J_{3,k}(t)$ and $\nabla J_{4,k}(t)$ are obtained in the same way.

The discretized (in space) set of optimality conditions result in a Two-Point Boundary Value Problem (TPBVP) for a coupled system of ODEs with time as independent variable. In order to fully solve the problem, we follow a standard iterative procedure to compute the reduced gradient~\cite{herzog}. We resort to a second-order time discretization to solve sequentially the state equation forward in time and the adjoint equation backward in time. Then, we evaluate Equation (\ref{OtD_grad}) to get the reduced gradient.

The fully discretized problem is obtained by applying the Crank-Nicolson method to the resulting system of ODEs.
We partition the time interval $[0,T]$ in $N$ sub-intervals of equal size $\Delta t = \frac{T}{N}$ and denote by $t_i = i\Delta t$ $i=0,\ldots,N$ the discrete time instances.  Furthermore, we denote the approximation of the unknown variable in the ODE system at time $t_i$ as $\mbf{q}(t_i) \approx \mbf{q}_i$. The Crank-Nicolson method for the state equation gives:                          

\begin{equation*}
\begin{aligned}
&\left[ \frac{M}{\Delta t}  +\frac{1}{2}\left\{ A + \Gamma_q(\ut_{i+1}) \right\}\right] \mathbf{q}_{i+1}   = \left[ \frac{M}{\Delta t}  -\frac{1}{2}\left\{ A +\Gamma_q(\ut_i)\right\} \right]\mathbf{q}_{i}   \\
& \bb{q}_0 = \bar{\bb{q}}_0,
\end{aligned}
\end{equation*}
for $i=0,\ldots,N-1$.
Defining the discrete transition matrices as:
\begin{equation}
\label{A_state}
\begin{aligned}
&\tilde{A}_{+}(\mathbf{u})  = \left[ \frac{M}{\Delta t}  +\frac{1}{2}\left\{ A + \Gamma_q(\ut) \right\}\right] \\
&\tilde{A}_{-}(\mathbf{u}) = \left[ \frac{M}{\Delta t}  -\frac{1}{2}\left\{ A +\Gamma_q(\ut)\right\} \right],
\end{aligned}
\end{equation}
we obtain the compact form:
\begin{equation}
\label{state_compact}
\begin{aligned}
    &\tilde{A}_{+}(\mathbf{u}_{i+1})\mathbf{q}_{i+1} = \tilde{A}_{-}(\mathbf{u}_{i})\mathbf{q}_{i} \quad , \quad i=0,\ldots,N-1 \\
    & \bb{q}_0 =\bar{\bb{q}}_0.
\end{aligned}
\end{equation}
In the same way, we define the discrete adjoint transition matrices
\begin{equation*}
\begin{aligned}
&\hat{A}_{+}(\mathbf{u}) = \left[ \frac{M}{\Delta t}   +\frac{1}{2}\left\{ A +\Gamma_p(\ut)\right\} \right] \\
&\hat{A}_{-}(\mathbf{u})  = \left[ \frac{M}{\Delta t}  -\frac{1}{2}\left\{ A + \Gamma_p(\ut) \right\}\right] 
\end{aligned}
\end{equation*}
so that the compact form of the adjoint dynamics reads as:
\begin{equation}
    \begin{cases}
    \label{adj_compact_OtD}
    &\hat{A}_{+}(\mathbf{u}_i) \mbf{p}_i = \hat{A}_{-}(\mathbf{u}_{i+1}) \mbf{p}_{i+1} \quad , \quad i=N-1,\ldots,0 \\
    & \mbf{p}_N = \mbf{q}_N - \bar{\mbf{q}}_T
    \end{cases}
\end{equation}
that has to be integrated backward in time.

\subsection{Discretize then Optimize}
The Discretize then Optimize (DtO) approach casts the OCP as a Nonlinear Programming Problem (NLP) where the dynamics is considered as a constraint. The set of constraints is already given by Equation (\ref{state_compact}) and we just need to fully discretize (i.e. both in time and space) the cost functional. Since the control functions are approximated in space with the RBF basis functions the cost functional reads as:
\begin{equation*}
\begin{aligned}
    J &= \frac{1}{2} \int_{\Omega}  (q(\x,T)-q_T(\x))^2 \, d \Omega    + \frac{\alpha}{2}\int_0^T \sum_{k=1}^4 \sum_{i=1}^{N_c} u_{k,i}^2(t) \, dt.
    \end{aligned}
\end{equation*}
Substituting the FEM approximation in space and approximating the integral in time with the trapezoidal method, that is consistent with the Crank-Nicolson method used to approximate in time the system of ODEs, we get:
\begin{equation*}
\begin{aligned}
    \tilde{J} &= \frac{1}{2} (\q_N-\bar{\q}_T)^{\top} M (\q_N-\bar{\q}_T)    \\ 
    &  +\frac{\alpha}{2}   \frac{\Delta t}{2}  \left( \sum_{i=1}^{N-1}2 \ut_i^{\top}\ut_i   + \ut_0^{\top}\ut_0+\ut_{N}^{\top}\ut_{N} \right);
    \end{aligned}
\end{equation*}
note that $\tilde{J} = \tilde{J}(\q_N,\ut_0,\ldots,\ut_N)$ is a \emph{function} of the control vector at each time instant and of the final state $\q_N$. This latter is also a function of the control vector through the fully discretized dynamics expressed in Equation (\ref{state_compact}). Therefore, the OCP can be cast as the following Nonlinear Program:
\begin{equation*}
\begin{aligned}
\min \quad & \tilde{J}(\q_N,\ut_0,\ldots,\ut_N)                 \\
\phantom{space}  \\
\textrm{s.t.} \quad &       \tilde{A}_{+}(\mathbf{u}_{i+1})\mathbf{q}_{i+1} = \tilde{A}_{-}(\mathbf{u}_{i})\mathbf{q}_{i} \quad i=0,\ldots,N-1 \\
& \bb{q}_0 = \bar{\bb{q}}_0  \\
& \mbf{0} \leq \ut_i \leq \ut_{max}   \quad i=0,\ldots,N.
\end{aligned}
\end{equation*}

We now use the discrete adjoint method \cite{schlog} to eliminate the dynamic constraints and to recast the NLP in the control unknowns $\ut_0,\ldots,\ut_N$ only. We define the \emph{discrete} Lagrangian \emph{function} as:

\begin{equation}
    \label{discrete_lag}
    \tilde{\mathcal{L}} = \tilde{J} - \Delta t \, \left\{ \sum_{i=0}^{N} \langle \tilde{A}_{+}(\ut_{i+1}) \q_{i+1} - \tilde{A}_{-}(\ut_{i}) \q_{i} , \mathbf{p}_{i+1} \rangle \right\}
\end{equation}
where $\bb{p}_i$ are the discrete adjoint variables associated to the dynamic constraint. Equation (\ref{discrete_lag}) can be rewritten as: 
\begin{equation}
    \label{discrete_lag_2}
\tilde{\mathcal{L}} = \tilde{J} - \Delta t \, \left\{ \sum_{i=1}^{N} \langle \tilde{A}_{+}(\ut_{i}) \q_{i} , \mathbf{p}_{i} \rangle - \sum_{i=0}^{N-1}\langle \tilde{A}_{-}(\ut_{i}) \q_{i} , \mathbf{p}_{i+1} \rangle \right\}.
\end{equation}
The gradient of the reduced cost functional $\tilde{J}(\q_N(\ut),\ut)$ is equal to the partial derivative of the discrete Lagrangian with respect to the control variables. This derivative for internal time instants (i.e. $i\neq 0,N$) is:
\begin{equation*}
\begin{split}
     &\frac{\partial\tilde{\mathcal{L}}}{\partial \mathbf{u}_i}^{\top} = \Delta t \,  \alpha \, \ut_i \\ 
     &-    \Delta t  \frac{\partial }{\partial \ut_i}\left\{  \langle \tilde{A}_{+}(\ut_{i}) \q_{i} , \mathbf{p}_{i} \rangle - \langle \tilde{A}_{-}(\ut_{i}) \q_{i} , \mathbf{p}_{i+1} \rangle \right\} .
\end{split}
\end{equation*}
In general, computing the derivative of $\frac{\partial }{\partial \ut_i}\tilde{A}_{+}(\ut_{i})$ might be an involved task, yielding a tensor. However, the matrix $\tilde{A}_{+}(\ut_{i}) $ is \emph{linear} in the control components and the derivative can be carried out term by term. Denoting $u_{i,1,k}$ the control component relative to the $k$-th basis function of actuator $1$ at time instant $i$ we obtain
\begin{equation*}
  \begin{aligned}
  &\frac{\partial \tilde{A}_{+}(\ut_{i})}{\partial u_{i,1,k}}  = \frac{\partial}{\partial u_{i,1,k}} \left[ \frac{M}{\Delta t}  +\frac{1}{2}\Big\{ A + \Gamma_q(\ut_i) \Big\}\right] \\
  &= \frac{\partial}{\partial u_{i,1,k}} \frac{1}{2} \Gamma_q(\ut_i) = \frac{1}{2} \Big\{B_{1,k}+C_{1,k}+L_{1,k}\Big\},
  \end{aligned}
\end{equation*}
then resulting in a sum of \emph{constant} matrices. The same results hold for the actuator sides $2,3,4$. The $\{ i,1,k\}$-th component of the reduced gradient of the cost function is:
\begin{equation*}
\begin{aligned}
    &\frac{\partial \tilde{J}}{\partial u_{i,1,k}} = \frac{\partial \tilde{\mathcal{L}}}{\partial u_{i,1,k}}   = \, \Delta t \, \Big(  \alpha \, u_{i,1,k}  \\
    &- \q_{i}^{\top}\, \Big( B^{\top}_{1,k}+C^{\top}_{1,k}+L^{\top}_{1,k} \Big) \frac{\left(\p_{i} + \p_{i+1}\right)}{2}  \, \Big)
    \end{aligned}
\end{equation*}
and, using again Lemma \ref{BCL}, we obtain:
\begin{equation}
\label{DtO_grad}
\begin{aligned}
    &\frac{\partial \tilde{J}}{\partial u_{i,1,k}} = \Delta t \, \Big(  \alpha \, u_{i,1,k} + \q_{i}^{\top}\, B_{1,k} \frac{\left(\p_{i} + \p_{i+1}\right)}{2}  \, \Big).
    \end{aligned}
\end{equation}
The components of the gradient for the control variables corresponding to the  actuator sides $2,3,4$ are obtained in a similar way. The components of the reduced gradient are readily obtained once we have the state and adjoint vectors $\q_i$ and $\mbf{p}_i$. The discrete adjoint equation is obtained by imposing that the derivative of the discrete Lagrangian with respect to the state vector $\q_i$ at each instant vanishes, that is:
\begin{equation*}
\begin{aligned}
        &\frac{\partial\tilde{\mathcal{L}}}{\partial \mathbf{q}_N}^{\top} = M  (\q_{N}-\q_T)
- \Delta t \, \tilde{A}_{+}^{\top}(\ut_{N}) \, \mathbf{p}_{N} = \mathbf{0},\\
&     \frac{\partial\tilde{\mathcal{L}}}{\partial \mathbf{q}_i}^{\top} = 
- \Delta t \left\{ \,  \tilde{A}_{+}^{\top}(\ut_{i})\mathbf{p}_{i} - \,  \tilde{A}^{\top}_{-}(\ut_{i})\mathbf{p}_{i+1} \right \} = \mathbf{0}. \\
\end{aligned}
\end{equation*}
The previous equations give the adjoint dynamics:
\begin{equation}
\label{adj_compact_DtO}
\begin{aligned}
 &\tilde{A}_{+}^{\top}(\ut_{i})\mathbf{p}_{i} =   \tilde{A}^{\top}_{-}(\ut_{i})\mathbf{p}_{i+1} \quad i=N-1,\ldots,1 \\
 &  \tilde{A}_{+}^{\top}(\ut_{N}) \, \mathbf{p}_{N} =  \frac{M}{\Delta t}  (\q_{N}-\bar{\q}_T) \\
\end{aligned}
\end{equation}
\subsection{Comparison between DtO and OtD approach}
DtO and OtD do not commute in general unless the dynamics is linear and the cost functional is quadratic~\cite{fredi,quart,schlog}. For the case at hand, the cost functional is quadratic, but the dynamics is bilinear since the control input multiplies the state. The bilinearity is present in the continuous setting in the state dynamics given by the combination of advection and reaction terms $\bb{v} \cdot \nabla q + (\nabla \cdot \bb{v}) q$. The FEM discretization is consistent with the PDE model and results in a bilinear system of ODEs since the matrix $\Gamma_q$ depends linearly on the control components. Since the problem is \emph{almost} linear we can show that the DtO and OtD \emph{almost} commute. We start drawing a comparison between the adjoint dynamics obtained with the two different approaches. The adjoint dynamics of the OtD approach is given by Equation (\ref{adj_compact_OtD}) while the DtO counterpart is given by Equation (\ref{adj_compact_DtO}).  Using Equation (\ref{gamma_qp}) it is easy to show that:
\begin{equation*}
    \tilde{A}_{+}^{\top}(\ut) = \hat{A}_{+}(\mathbf{u}).
\end{equation*}
Indeed, from its definition in Equation (\ref{A_state}) we have: 
\begin{equation*}
\begin{aligned}
&\tilde{A}^{\top}_{+}(\mathbf{u})  = \left[ \frac{M^{\top}}{\Delta t}  +\frac{1}{2}\left\{ A^{\top} + \Gamma_q(\ut)^{\top} \right\}\right]  \\
&=\left[ \frac{M}{\Delta t}  +\frac{1}{2}\left\{ A+ \Gamma_q(\ut)^{\top} \right\}\right]  \\
&=\left[ \frac{M}{\Delta t}  +\frac{1}{2}\left\{ A+ \Gamma_p(\ut) \right\}\right] = \hat{A}_{+}(\mathbf{u}).
\end{aligned}
\end{equation*}
Hence, the two different formulations of the adjoint dynamics in Equations (\ref{adj_compact_OtD}) and (\ref{adj_compact_DtO}) differ only at the time instants at which they are evaluated. The reduced gradients obtained with the two approaches in Equations (\ref{OtD_grad}) and (\ref{DtO_grad}) only differ for the time instant which the adjoint is evaluated at. In the OtD formulation, everything is continuous in time and each variable is evaluated at the same instant when the time is discretized. In the DtO approach, the adjoint is evaluated at the midpoint between two successive time instants, while state and control functions are evaluated at the same time instant.

The main result that allows to establish the \emph{quasi} commutation of OtD and DtO approaches is Lemma \ref{BCL} that allows to prove Equation \eqref{gamma_qp}).  Thanks to this it is possible to infer that the semi-discretization in space fully commutes. The semi-discrete adjoint equation arising from the OtD method is Equation \eqref{otd_adjoint} while it is possible to recover the semi-discrete adjoint equation arising from the DtO method with a passage to the limit. Starting from the adjoint algebraic Equation \eqref{adj_compact_DtO} and using the definition of $\tilde{A}_{+}^{\top}(\ut)$ we obtain:
\begin{equation*}
\begin{aligned}
&\tilde{A}^{\top}_{-}(\ut_{i})\mathbf{p}_{i+1} - \tilde{A}_{+}^{\top}(\ut_{i})\mathbf{p}_{i} \\
&= M \frac{\mathbf{p}_{i+1} - \mathbf{p}_{i}}{\Delta t} - \frac{1}{2}\big(A+\Gamma_q^{\top}(\mathbf{u}_i)\big) \mathbf{p}_{i+1} - \frac{1}{2}\big(A+\Gamma_q^{\top}(\mathbf{u}_i)\big) \mathbf{p}_{i} 
\end{aligned}
\end{equation*}
in the limit $\Delta t \to 0$ we formally have $ \frac{\mathbf{p}_{i+1} - \mathbf{p}_{i}}{\Delta t} = \dot{\mathbf{p}}$ and $\mathbf{p}_{i+1} \approx \mathbf{p}_{i} = \mathbf{p}$ and thence the semi-discrete adjoint dynamics resulting from the DtO method is:
\begin{equation*}
- M\dot{\mathbf{p}} + A \mathbf{p} + \Gamma^{\top}_q(\mathbf{u})\mathbf{p} = \mathbf{0}.
\end{equation*}
Therefore the two approaches fully commute (in space) due to Equation \eqref{gamma_qp} whose proof relies on Lemma \ref{BCL}.

\section{Numerical simulation}
\label{num_sim}

In this Section, the OCP is fully discretized in time using the Crank-Nicolson method and the resulting Nonlinear Optimization Program (NLP) is stated. The discrete adjoint method is then used to compute the gradient of the reduced objective function and to set up the optimization in the control space only.

The control coefficients at each instant of time and for each basis function are stacked in a single vector $\mathbf{U} = [\ut_0,\ldots,\ut_N] \in \R^{N_c (N+1)}$. In the same way. we define the stacked state and adjoint vector $\mathbf{Q} = [\bb{q}_0,\ldots,\bb{q}_N] \in \R^{N_q (N+1)}$ and $\mathbf{P} = [\bb{p}_1,\ldots,\bb{q}_N] \in \R^{N_p N}$.

The resulting optimization problem  has the form:

\begin{equation}
\label{final_opt}
\begin{aligned}
\min_{\mathbf{U}} \quad & \tilde{J}(\mathbf{U})                 \\
\textrm{s.t.} \quad &  \mathbf{0} \leq \mathbf{U} \leq \mathbf{U}_{max}.
\end{aligned}
\end{equation}

Problem (\ref{final_opt}) is a NLP subject to bound constraints only. The procedure to numerically compute the exact gradient $\nabla \tilde{J}(\mathbf{U})$ is given in Algorithm \ref{alg}.

\begin{algorithm}
  \begin{algorithmic}[1]
    \State Given $\mathbf{U} = [\ut_0,\ldots,\ut_N]$
    \State Solve state Equation (\ref{state_compact}) for  $\mathbf{Q}= [\mbf{q}_0,\ldots,\mbf{q}_N]$  
    \State Solve adjoint Equation (\ref{adj_compact_DtO}) for $\mathbf{P}= [\mbf{p}_1,\ldots,\mbf{p}_N]$
    \State Evaluate Equations (\ref{DtO_grad}) for  $\nabla \tilde{J} (\mathbf{U})$
  \end{algorithmic}
  \caption{Reduced Gradient  }
  \label{alg}
\end{algorithm}

We make use of a MATLAB interface \cite{mexipopt} to the NLP solver \texttt{IPOPT} \cite{ipopt} to solve the NLP (\ref{final_opt}). We provide the solver with the reduced gradient obtained from Algorithm (\ref{alg}).

In the next subsections we will show the numerical results obtained for two different test cases. For both of them the main parameters of the simulation are summarized in Table \ref{param_table}, so that the time instants are $N=40$, while, selecting $N_c=10$ basis functions, the stacked control vector $\mathbf{U}$ has dimension $1640$ while the state $\mathbf{Q}$ and adjoint vectors $\mathbf{P}$ have $29192$ and $28480$ elements respectively.
All computations were conducted on a Dell XPS15 desktop PC with an Intel Core i7-10750HQ CPU and 16 GB RAM running Ubuntu 18.04. The computational time is roughly 3 hours for both cases.

\subsection{Test case 1}
\begin{figure}[hbt!]
\centering
\includegraphics[width=3in]{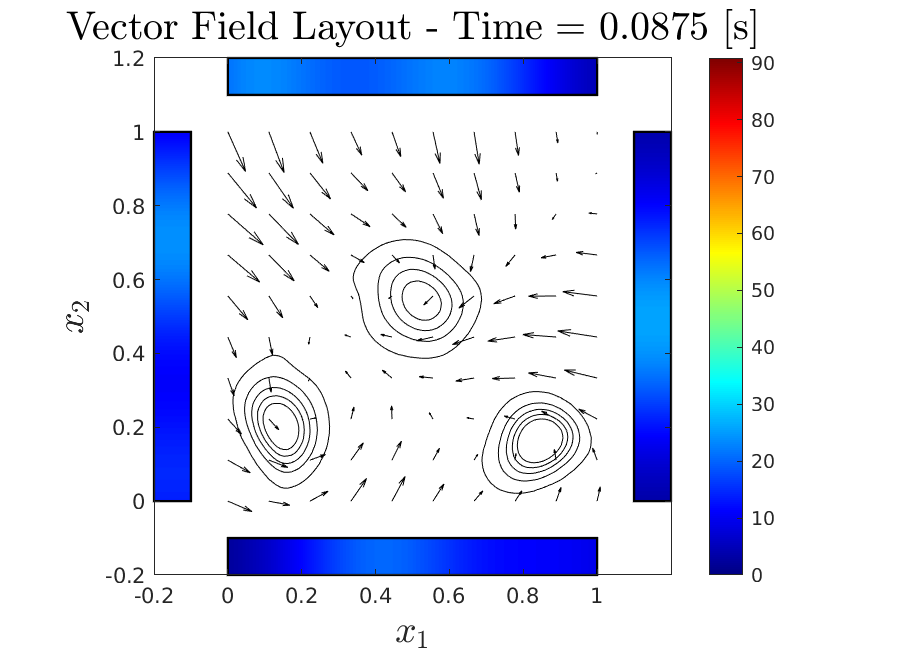}
\caption{Layout of the control system at $t=0.0875 \, [s]$. The actuator stacks spatial intensity is shown according to the colormap on the right. The induced velocity field is represented by black arrows in the workspace together with the level-set curves of the density. } 
\label{layout_colored}
\end{figure}

As a first test case, we consider the problem of driving an initial uniform density to a combination of three disjoint radial basis functions. Figure \ref{layout_colored} presents a snapshot of the control system where the spatial intensity of the actuator stacks is shown together with the induced vector field.

\begin{table}[h]
\begin{center}
\begin{tabular}{ |c|c|} 

 \hline
 Parameter & Value \\
 \hline 
 $T$ &    0.1    \\ 
 $\Delta t$     &    0.0025     \\ 
 $\alpha$      &     0.0001    \\
 $\mu$        &      0.1         \\
 FEM nodes & 712 \\
 \hline
\end{tabular}
\caption{Main simulation parameters. \label{param_table}} 
\end{center}
\end{table}

\begin{figure}[hbt!]
\centering
\includegraphics[width=\linewidth]{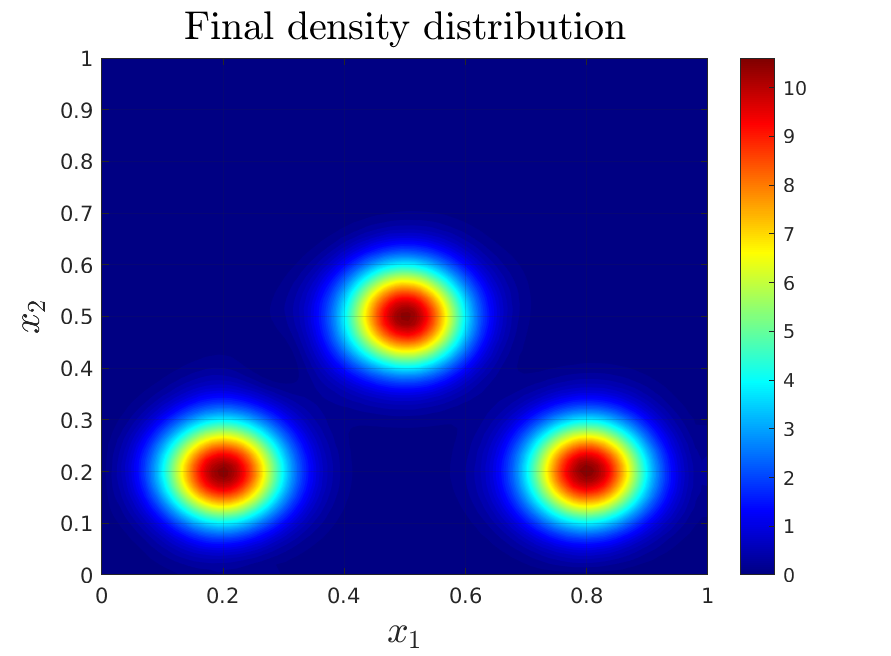}
\caption{Test Case 1. Density reached at final time $T=0.1 \, [s]$ by the optimal control algorithm.  } 
\label{final_density}
\end{figure}

\begin{figure}[hbt!]
\centering
\includegraphics[width=\linewidth]{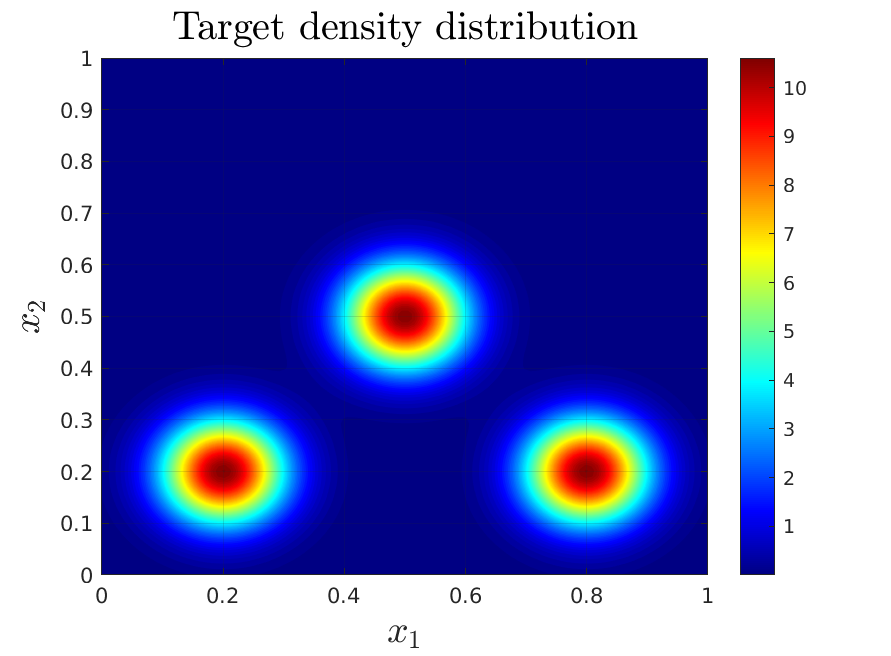}
\caption{Test Case 1. Target Density consisting of three normalized disjoint radial basis functions.} 
\label{target_density}
\end{figure}

\noindent The target density consists of a combination of Radial Basis Functions with unitary total mass so that they represent a probability density. The contour plot of the target density is shown in Figure \ref{target_density}. The density at the final instant from the solution of the OCP is shown in Figure \ref{final_density}.

\begin{figure}[hbt!]
\centering
\includegraphics[width=\linewidth]{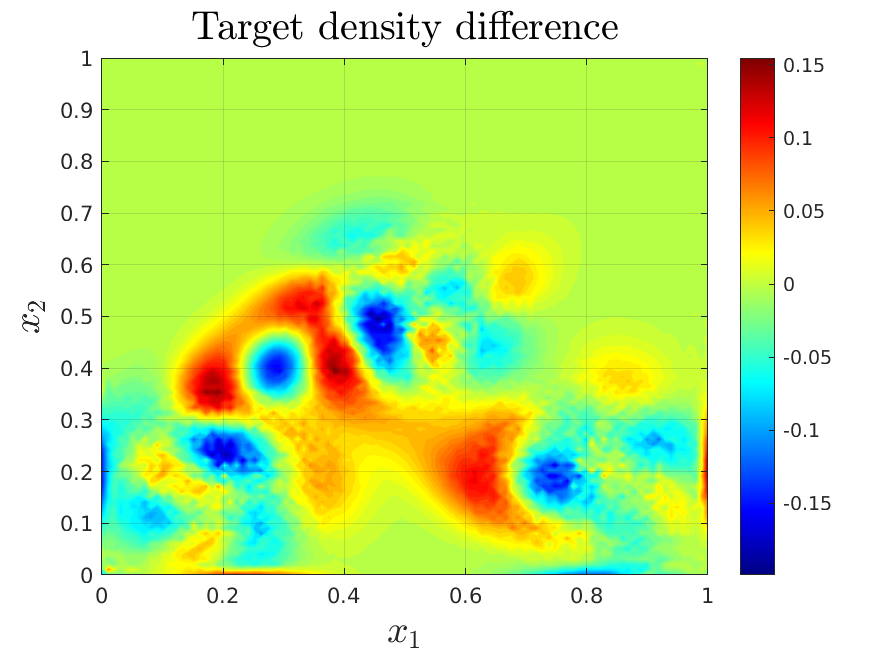}
\caption{Test Case 1. Pointwise difference between target density $q_T(\bb{x})$ and the optimal density $q^{\star}(\bb{x},T)$ reached at the final instant by the optimal control algorithm. The $L^2$ distance between $q^{\star}(\bb{x},T)$ and $q_T(\bb{x})$ is $\int_{\Omega} (q^{\star}(\bb{x},T)-q_T(\bb{x}))^2 \, d\Omega = 0.001$. } 
\label{target_diff}
\end{figure}

Figure \ref{target_diff} shows the difference between the target and the density reached at the final instant demonstrating that the optimal control algorithm is able to reach the target.

\begin{figure}[hbt!]
\centering
\includegraphics[width=\linewidth]{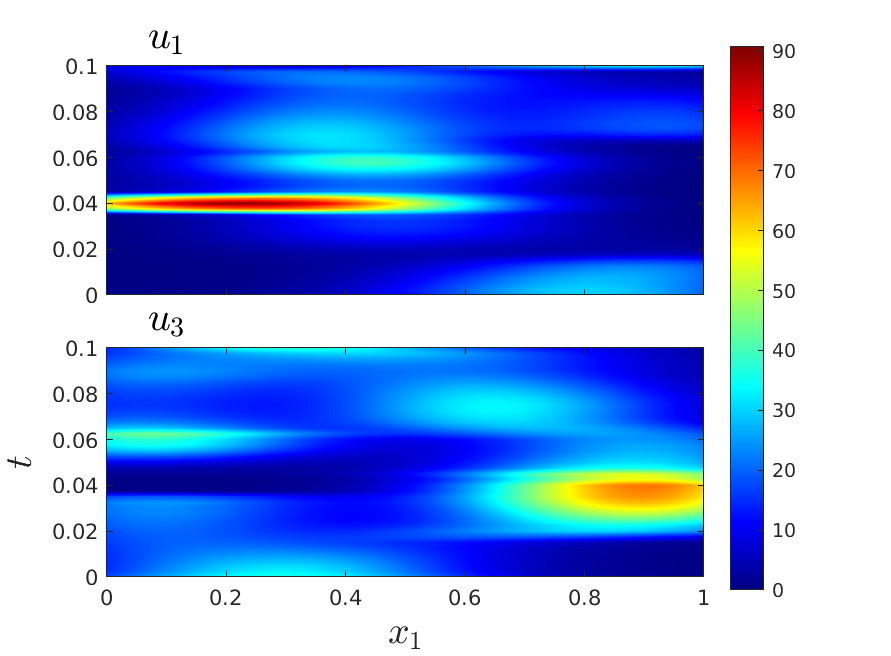}
\caption{Test Case 1. Space-time intensity of actuators $u_1$ and $u_3$ that generates a velocity field along the $x_2$ axis in the positive and negative direction respectively. } 
\label{u13_1}
\end{figure}

\begin{figure}[hbt!]
\centering
\includegraphics[width=\linewidth]{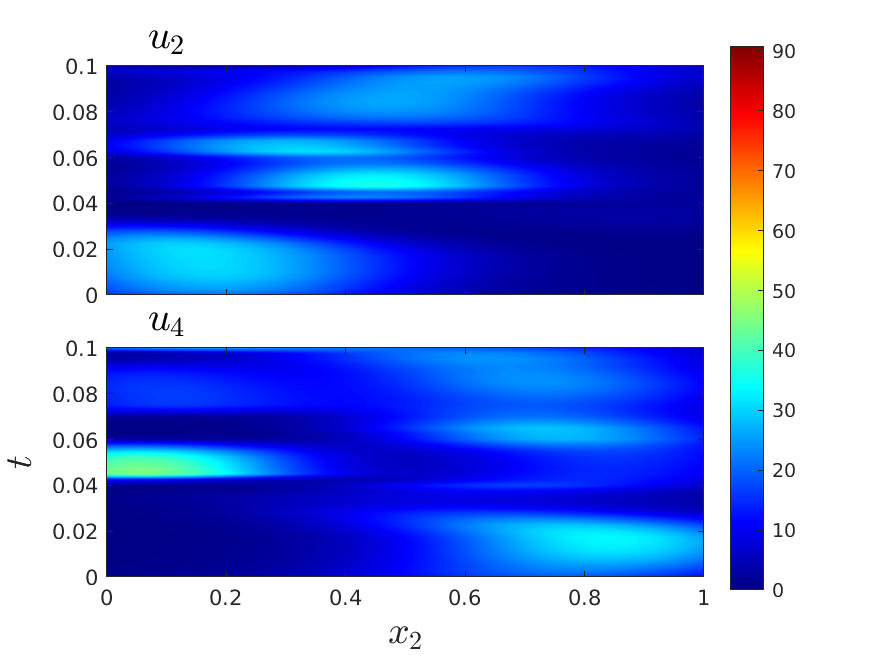}
\caption{Test Case 1. Space-time intensity of actuators $u_2$ and $u_4$ that generates a velocity field along the $x_1$ axis in the nevative and positive direction respectively. } 
\label{u24_1}
\end{figure}

Finally, the optimal control space-time history of the actuation in the vertical and horizontal direction is shown in Figure \ref{u13_1} and Figure \ref{u24_1}, respectively. It is interesting to notice that the maximum actuation effort is needed from the left portion of actuator $u_1$. Most of the space-time domain intensity of the actuators is close to zero due to the minimum energy feature of the optimal control problem.

\subsection{Test case 2}
\begin{figure}[hbt!]
\centering
\includegraphics[width=\linewidth]{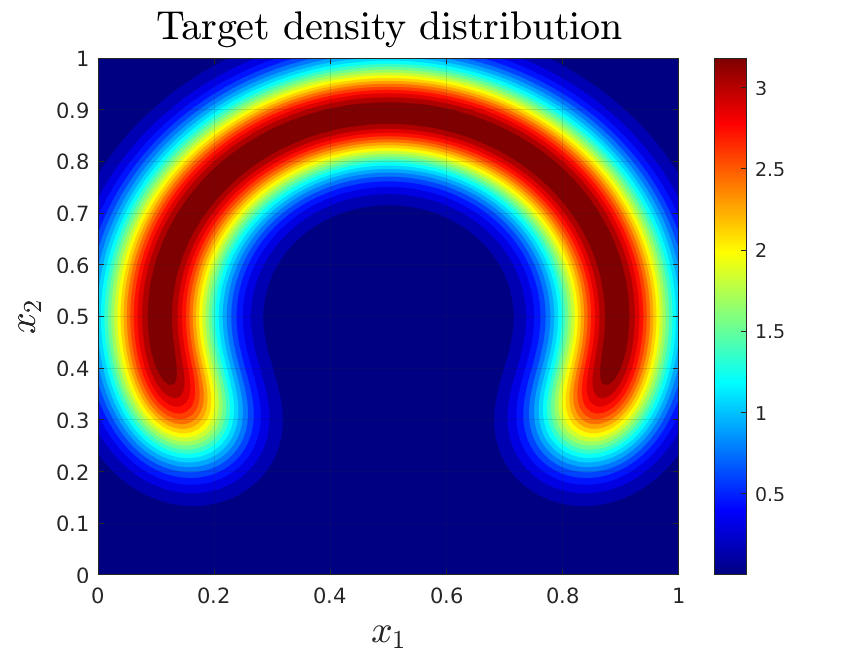}
\caption{Test Case 2. Target Density consisting of a complex horseshoe-shaped function.} 
\label{target_density_2}
\end{figure}

\begin{figure}[hbt!]
\centering
\includegraphics[width=\linewidth]{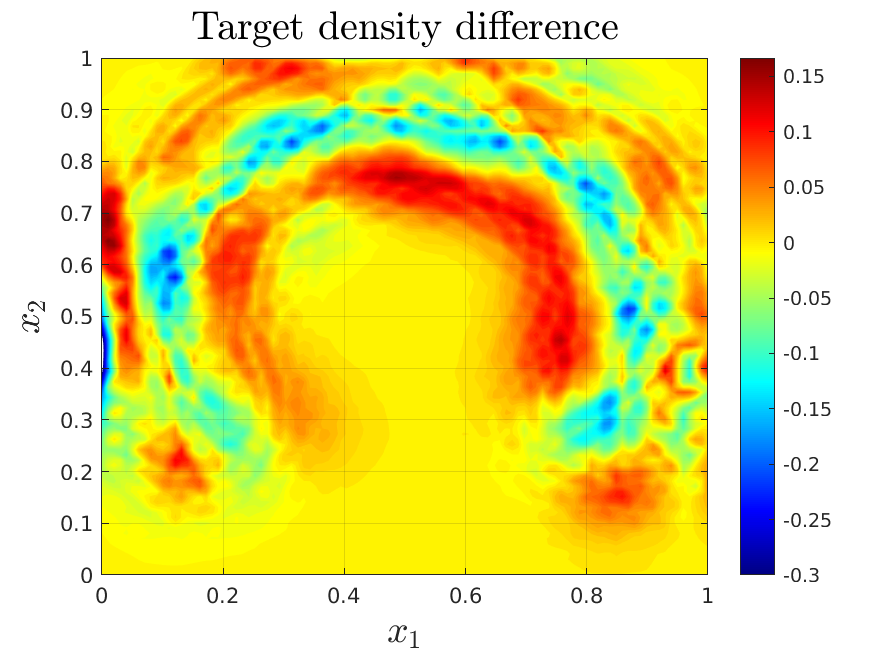}
\caption{Test Case 2. Pointwise difference between target density $q_T(\bb{x})$ and the optimal density $q^{\star}(\bb{x},T)$ reached at the final instant by the optimal control algorithm. The $L^2$ distance between $q^{\star}(\bb{x},T)$ and $q_T(\bb{x})$ is $\int_{\Omega} (q^{\star}(\bb{x},T)-q_T(\bb{x}))^2 \, d\Omega = 0.0015$. } 
\label{difference_2}
\end{figure}

\begin{figure}[hbt!]
\centering
\includegraphics[width=\linewidth]{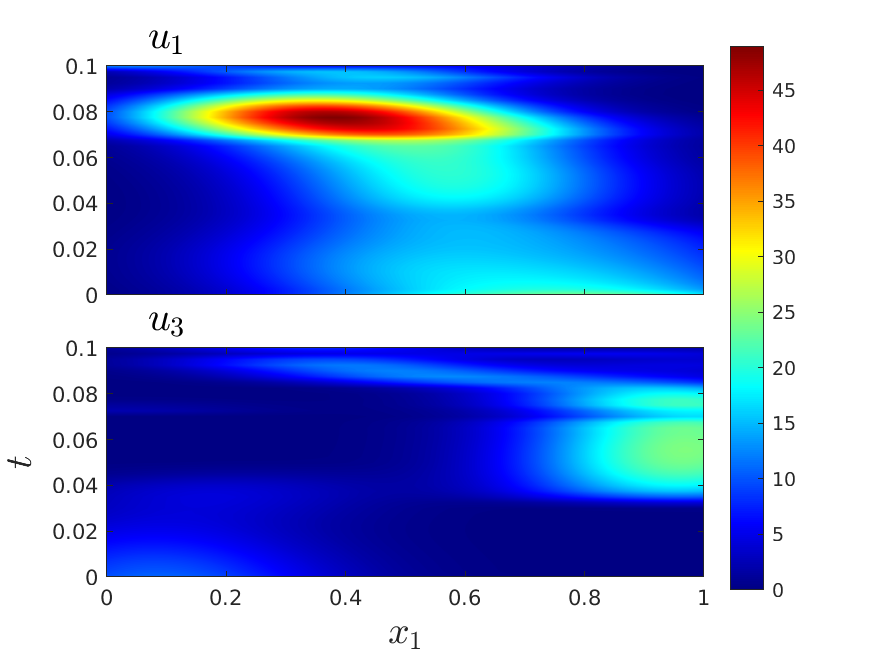}
\caption{Test Case 2. Space-time intensity of actuators $u_1$ and $u_3$.}
\label{u13_2}
\end{figure}

\begin{figure}[hbt!]
\centering
\includegraphics[width=\linewidth]{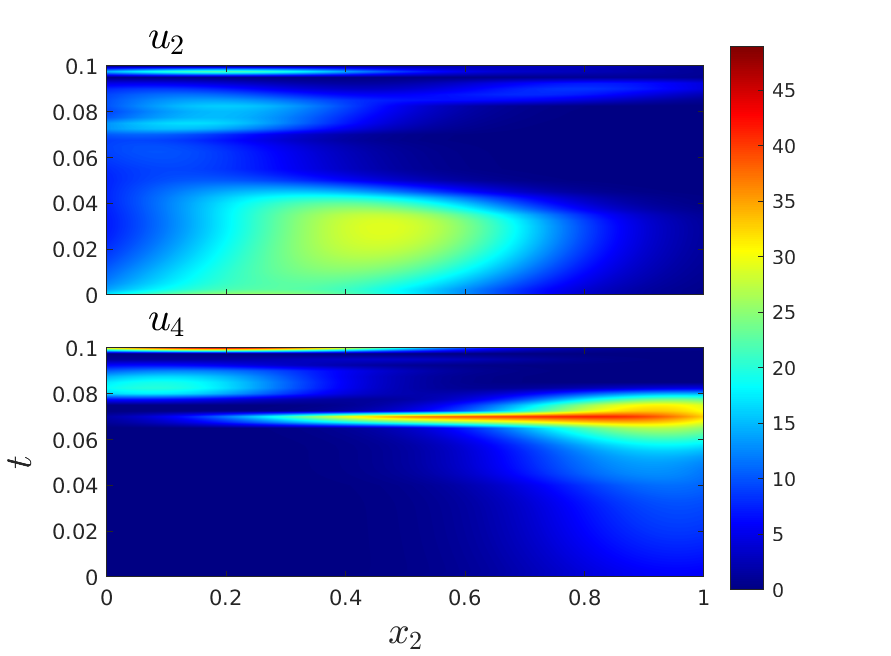}
\caption{Test Case 2. Space-time intensity of actuators $u_2$ and $u_4$.} 
\label{u24_2}
\end{figure}

\begin{figure}[hbt!]
\centering
\includegraphics[width=\linewidth]{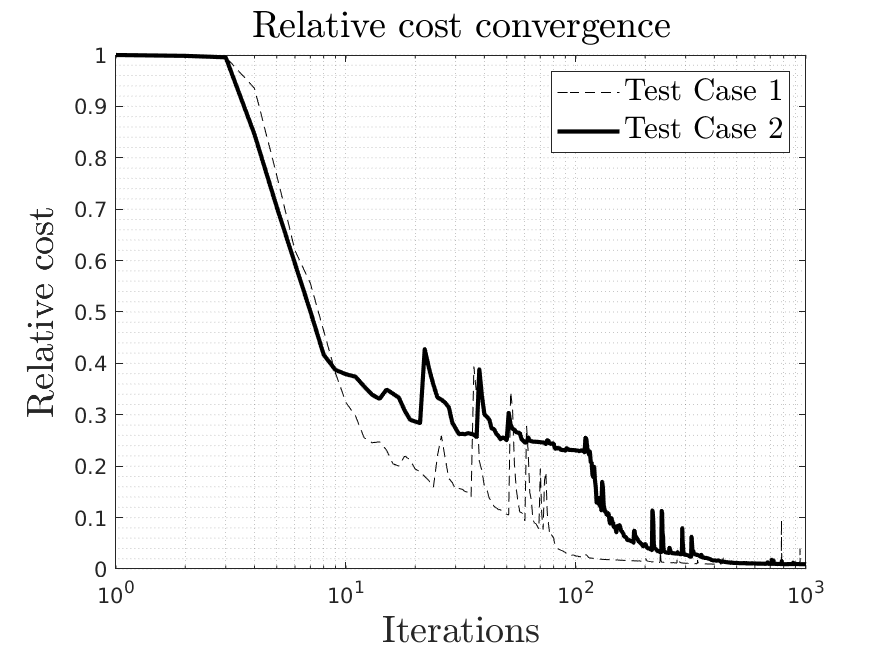}
\caption{Iterations of the NLP solver for both test cases. Each cost is normalized by the initial value of the cost functional and the abscissa axis is in logarithmic scale. For Test Case 1 we have $J^{0}=2.14$, while for Test Case 2 we have $J^{0}=0.66$.} 
\label{cost_comparison}
\end{figure}

As a second test case, we consider the problem of driving an initial uniform density to a complex target function shown in Figure \ref{target_density_2}. The control algorithm is able to drive the density to the target as shown in Figure \ref{difference_2}. The simulation parameters are the same as for Test Case $1$. We note that the control effort is significantly lower as shown in Figures \ref{u13_2} and \ref{u24_2}. This is due to the fact that the target density, despite being more complex, it is more evenly distributed in the workspace and it is closer in the $L^2$ sense to the initial uniform distribution. Finally, we report in Figure \ref{cost_comparison} a comparison between relative cost iterations. The initial guess for the control $\mathbf{U}^{0}$ is the zero vector. Note that the initial condition $q_0(\x)$ is a uniform density and it is an equilibrium distribution for the zero control case. Therefore, the value of the cost functional at the first iteration $J^{0}$ is the FEM approximation of $\frac{1}{2} \int_{\Omega} (q_T(\x)-q_0(\x))^2 \, d\Omega $ that is half the $L^2$ distance between $q_T$ and $q_0$. The relative cost iterations are slightly faster to converge for Test Case 1 thus showing that Test Case 2 represents a more difficult problem. The initial cost is decreased of $99.4\%$ and $99.1\%$ respectively.

\section{Conclusions}
\label{conclusion}
In this paper we have showed how to contain and control a large scale swarm of underactuated particles. Actually, the method presented assumes an infinite number of passive particles. The macroscopic dynamic model encodes the physical layout of the actuators and automatically takes into account the limited control authority. The PDE model arising from a Kolmogorov forward equation is used as state dynamics to set up an optimal control problem. The state dynamics and the optimal control problem are thoroughly analysed. A series of estimates is provided for the state dynamics that allowed to prove an existence theorem for the optimal control vector.  The necessary conditions for optimality for this problem are analytically derived in closed-form expression. From the continuous formulation, the two main numerical approaches (i.e. DtO and OtD) are investigated and the resulting discrete equations are compared. By exploiting some properties of the resulting matrices, it is shown that the two approaches almost commute. In addition, we note that the numerical procedure is general in nature and can work with any kind of control basis functions. Compared to previous works \cite{nc_fod,gen_grad}, our method allows the control basis function to be not necessarily null at the boundary thus improving the flexibility of the resulting systems. The adjoint sensitivity analysis used to compute the gradient of the reduced cost functional in the DtO approach allows for an exact and fast way to compute the sensitivity. Compared to similar results in the literature \cite{schlog}, we used the more accurate Crank-Nicolson method for the time discretization, carefully taking into account the state transition matrices resulting from the fully discretized system. The numerical simulation showed the effectiveness of the resulting method. Future research will mainly focus on two aspects. On the one hand, the link between macroscopic and microscopic dynamics will be further investigated considering a finite number of particles subjected to the optimal velocity field obtained. On the other hand, a mechanism to encode more complex particles behaviour such as repulsion will also be considered.

\appendix
\subsection{Derivation of the bilinear form \eqref{bil_form} }
The bilinear form $a(q,\phi)$ is defined as
\begin{equation*}
a(q,\phi) = \int_{\Omega}   \nabla \cdot (-\mu\nabla q + \bb{v} q) \phi  d\Omega,
\end{equation*}
using integration by parts and the boundary conditions we can show that:
\begin{equation*}
 \int_{\Omega}   \nabla \cdot (-\mu\nabla q + \bb{v} q) \phi  d\Omega =  \int_{\Omega} \mu \nabla q \cdot \nabla \phi - \bb{v} \, q \cdot \nabla \phi d\Omega.
\end{equation*}
Indeed:
\begin{equation*}
\begin{aligned}
&\int_{\Omega}   \nabla \cdot (-\mu\nabla q + \bb{v} q) \phi \,d\Omega  = \int_{\Omega} \mu \nabla q \cdot \nabla \phi \,d\Omega \\
& - \int_{\Gamma} \mu \nabla q \cdot \nor \, \phi \, d\Omega + \int_{\Omega} \nabla \cdot (\bb{v} q) \phi \,d\Omega  =  \int_{\Omega} \mu \nabla q \cdot \nabla \phi\,d\Omega \\
& - \int_{\Gamma} \mu \nabla q \cdot \nor \, \phi\,d\Omega  + \int_{\Gamma} \bb{v} \cdot \nor q \phi\,d\Omega \\
& - \int_{\Omega} q \bb{v} \cdot \nabla \phi \,d\Omega  = \int_{\Omega}   \mu \nabla q \cdot \nabla \phi - q\,\bb{v} \cdot \nabla \phi\,d\Omega \\
&+ \cancel{\int_{\Gamma} (-\mu \nabla q + \bb{v} q ) \cdot \nor  \, \phi \,d\Gamma }
\end{aligned}
\end{equation*}
so that the bilinear form can be equivalently written as:
\begin{equation*}
a(q,\phi) =  \int_{\Omega} \mu \nabla q \cdot \nabla \phi - \bb{v} \, q \cdot \nabla \phi.
\end{equation*}

\subsection{Proof of Lemma \ref{lemma_1}}
\label{appendix_v}
We define:
\begin{equation*}
\begin{split}
    &\mathbf{b}_1 = \begin{bmatrix}
    0 \\
    e^{-cx_2}
    \end{bmatrix} ,
\quad
    \mathbf{b}_2 = \begin{bmatrix}
    0 \\
    e^{-c(1-x_1)}
    \end{bmatrix} ,
    \\
    &\mathbf{b}_3 = \begin{bmatrix}
    0 \\
    e^{-c(1-x_2)}
    \end{bmatrix} , 
\quad
    \mathbf{b}_4 = \begin{bmatrix}
    0 \\
    e^{-cx_1}
    \end{bmatrix} ;
\end{split}
    \end{equation*}
hence, we can express $\bb{v}(t)$ as:
\begin{equation*}
    \bb{v}(t) = \sum_{i=1}^4 \bb{b}_i u_i(t).
\end{equation*}
For every $t>0$, we thus have:
\begin{equation*}
\begin{split}
&\norm{\bb{v}(t)}_{L^{\infty}(\Omega)^2} = \norm{\sum_{i=1}^4 \bb{b}_i u_i(t) }_{L^{\infty}(\Omega)^2} \\ 
&\leq \sum_{i=1}^4 \norm{\bb{b}_i u_i(t) }_{L^{\infty}(\Omega)^2} \leq \sum_{i=1}^4 \norm{ u_i(t) }_{L^{\infty}(\Gamma_i)}
\end{split}
\end{equation*}
Squaring and integrating in time between $0$ and $T$; we obtain
\begin{equation*}
\begin{split}
    &\norm{\bb{v}}^2_{L^2(0,T;L^{\infty}(\Omega)^2)}  = \int_0^T \norm{\bb{v}(t)}^2_{L^{\infty}(\Omega)^2} \\ &\leq \int_0^T \Big(\sum_{i=1}^4 \norm{ u_i(t) }_{L^{\infty}(\Gamma_i)}\Big)^2 \\
    &\leq 8 \sum_{i=1}^4 \int_0^T   \norm{ u_i(t) }^2_{L^{\infty}(\Gamma)} =  8 \norm{\mathbf{u}}^2_{\mathcal{U}},
\end{split}
\end{equation*}
where we iteratively used Cauchy inequality $(a+b)^2 \leq 2 (a^2+b^2)$. Then, turning to the definition of norms we finally obtain:
\begin{equation*}
    \norm{\bb{v}}^2_{L^2(0,T;L^{\infty}(\Omega)^2)} \leq 8 \norm{\mathbf{u}}^2_{\mathcal{U}}.
\end{equation*}
\QEDA

\subsection{Proof of Lemma \ref{weakly_coerc}}
Using the definition of $a(q,\phi)$ and Cauchy-Schwarz inequality we have
\begin{equation*}
\begin{aligned}
& a(q(t),q(t)) +\lambda(t) \norm{q(t)}^2_{L^2(\Omega)} \geq \mu \norm{\nabla q(t)}_{L^2(\Omega)}^2   \\
& - \norm{\bb{v}(t)}_{L^{\infty}(\Omega)^2} \norm{\nabla q(t)}_{L^2(\Omega)}  \norm{ q(t) }_{L^2(\Omega)} +\lambda(t) \norm{q(t)}^2_{L^2(\Omega)} \\ &\geq  (\mu-\epsilon(t))\norm{\nabla q(t)}_{L^2(\Omega)}^2\\
&+ \Big(\lambda(t)-\frac{\norm{\mathbf{v}(t)}^2_{L^{\infty}(\Omega)^2}}{4\, \epsilon(t)}\Big)\norm{ q(t) }^2_{L^2(\Omega)}.
\end{aligned}
\end{equation*}
Note that we have used Cauchy's inequality (see e.g., \cite{evans}), that is $ab \leq \epsilon a^2 + \frac{b^2}{4\epsilon}$ for $a,b,\epsilon > 0$ so that we have:
\begin{equation*}
\begin{aligned}
&\norm{ \nabla q(t)}_{L^2(\Omega)}\big(\norm{\bb{v}(t)}_{L^{\infty}(\Omega)^2} \, \norm{q(t)}_{L^2(\Omega)}  \big)  \\
& \leq  \epsilon \norm{ \nabla q(t)}_{L^2(\Omega)}^2 + \frac{\norm{\bb{v}(t)}_{L^{\infty}(\Omega)^2}^2 \, \norm{q(t)}_{L^2(\Omega)}^2 }{4 \epsilon}
\end{aligned}
\end{equation*}

Then, setting $\epsilon(t) = \frac{\mu}{2}$, we need to choose $\lambda(t)$ such that $\lambda(t) > \frac{\norm{\bb{v}(t)}^2_{L^{\infty}(\Omega)^2}}{2 \mu} $; hence, it is sufficient to select
\begin{equation*}
    \lambda(t) = \frac{\norm{\mathbf{v}(t)}^2_{L^{\infty}(\Omega)^2}}{\mu}.
\end{equation*}
Finally choosing  $\alpha_0(t) = \min\left\{\frac{\mu}{2},\frac{\norm{\bb{v}(t)}^2_{L^{\infty}(\Omega)^2}}{2 \mu} \right\}$ we have
\begin{equation*}
\begin{split}
    & a(q(t),q(t)) +\lambda(t) \norm{q(t)}^2_{L^2(\Omega)}  \\ & \geq \frac{\mu}{2}\norm{\nabla q(t)}_{L^2(\Omega)}^2 + \frac{\norm{\mathbf{v}(t)}^2_{L^{\infty}(\Omega)^2}}{2\mu}\norm{ q(t) }_{L^2(\Omega)}\\ 
    & \geq  \alpha_0(t) \norm{q(t)}_{H^1(\Omega)}^2.
\end{split}
\end{equation*}
\QEDA

We are now ready to prove Theorem \ref{theorem_1}.
\subsection{Proof of Theorem \ref{theorem_1}}

\subsubsection{Existence and uniqueness}
From Lemma \ref{lemma_1} we have that $\bb{v} \in L^{2}(0,T;L^{\infty}(\Omega)^2)$; from Lemma \ref{weakly_coerc}, the bilinear form $a(q,\phi)$ associated to the weak formulation of the problem (\ref{weak_form}) is weakly coercive. 
$a(q(t),\phi(t))$ is continuous in $H^1(\Omega)$ for a.e. $t \in (0,T)$ since
\begin{equation*}
|a(q(t),\phi(t))| \leq \Big( \mu + \norm{\bb{v}(t)}_{L^{\infty}(\Omega)^2} \Big) \norm{q(t)}_{H^1(\Omega)} \norm{\phi(t)}_{H^1(\Omega)}
\end{equation*}
and t-measurable for fixed $q(t)$,$\phi(t)$ and $\bb{v}(t)$. That is $a: t \mapsto a(q(t),\phi(t))$ is in $L^1(0,T)$. These assumptions guarantee the well-posedness of the state problem and thus the existence of a unique $q \in L^2(0,T;H^1(\Omega))$ which solves  (\ref{state_eq}) with $\dot{q} \in L^2(0,T;H^{1}(\Omega)^*)$, see, e.g., Chapter 7 in \cite{MQS}.

\subsubsection{Proof of Estimate (\ref{subeq1}) and (\ref{subeq2})}
Substituting $q(t)$ as test function in the weak formulation in Equation (\ref{weak_form}) we have:
\begin{equation*}
         \int_{\Omega} \frac{\partial q(t)}{\partial t} q(t) d \Omega + a(q(t),q(t)) = 0 \quad a.e. \,\, t \in (0,T),
\end{equation*}
that can be rewritten equivalently as
\begin{equation}
        \label{weak_form_2}
         \frac{1}{2}\frac{d}{dt}\norm{q(t)}^2_{L^2(\Omega)} + a(q(t),q(t)) = 0 \quad a.e. \,\, t \in (0,T).
\end{equation}

From Lemma \ref{weakly_coerc}, Equation (\ref{weak_form_2})  gives:
\begin{equation*}
\begin{split}
    &\frac{1}{2}\frac{d}{dt}\norm{q(t)}^2_{L^2(\Omega)} + a(q(t),q(t)) \\
    & \geq \frac{1}{2}\frac{d}{dt}\norm{q(t)}^2_{L^2(\Omega)} + \alpha_0 \norm{q(t)}^2_{H^1(\Omega)} - \lambda(t) \norm{q(t)}^2_{L^2(\Omega)} 
\end{split}
\end{equation*}
which implies:
\begin{equation*}
\frac{1}{2}\frac{d}{dt}\norm{q(t)}^2_{L^2(\Omega)}+\alpha_0(t) \norm{q(t)}^2_{H^1(\Omega)} \leq \lambda(t) \norm{q(t)}^2_{L^2(\Omega)}   .
\end{equation*}
Integrating between $0$ and $t$; we obtain
\begin{equation}
\begin{split}
    \label{int_norm}
    &\frac{1}{2}\norm{q(t)}^2_{L^2(\Omega)}-\frac{1}{2}\norm{q_0}^2_{L^2(\Omega)}+\int_0^t \alpha_0(\tau) \norm{q(\tau)}^2_{H^1(\Omega)}\,d\tau \\
    &\leq \int_0^t \lambda(\tau)  \norm{q(\tau)}^2_{L^2(\Omega)} \, d\tau,
\end{split}
\end{equation}
so that
\begin{equation*}
    \norm{q(t)}^2_{L^2(\Omega)} \leq \norm{q_0}^2_{L^2(\Omega)} +\int_0^t 2\,\lambda(\tau)  \norm{q(\tau)}^2_{L^2(\Omega)} \, d\tau.
\end{equation*}
Applying Gronwall Lemma (see, e.g., \cite{evans}) we have
\begin{equation}
\label{gron_1}
\norm{q(t)}^2_{L^2(\Omega)}  \leq e^{\int_0^t 2\,\lambda(\tau)  d\tau} \norm{q_0}^2_{L^2(\Omega)}.
\end{equation}
We now look for a bound on the term $\int_0^t 2\,\lambda(\tau)  d\tau$. From the proof of Lemma \ref{lemma_1}, we have that:
\begin{equation*}
    \norm{\bb{v}(t)}^2_{L^{\infty}(\Omega)^2} \leq  8 \sum_{i=1}^4 \norm{ u_i(t) }^2_{L^{\infty}(\Gamma_i)}
\end{equation*}
recalling the expression of $\lambda(t)$ found in Lemma \ref{weakly_coerc} we have:
\begin{equation*}
    \lambda(t) = \frac{\norm{\mathbf{v}(t)}^2_{L^{\infty}(\Omega)^2}}{\mu} \leq  \frac{8}{\mu} \sum_{i=1}^4 \norm{ u_i(t) }^2_{L^{\infty}(\Gamma_i)} 
\end{equation*}
and hence
\begin{equation*}
    \int_0^t 2\,\lambda(\tau)  d\tau \leq \int_0^T 2\,\frac{8}{\mu} \sum_{i=1}^4 \norm{ u_i(\tau) }^2_{L^{\infty}(\Gamma_i)} \, d\tau = \frac{16}{\mu} \norm{\mathbf{u}}^2_{\mathcal{U}}
\end{equation*}
We can then bound the terms in Equation (\ref{gron_1}) from above as:
\begin{equation}
    \label{base_est}
    \norm{q(t)}_{L^2(\Omega)}^2  \leq e^{\frac{16}{\mu} \norm{\mathbf{u}}^2_{\mathcal{U}}} \norm{q_0}_{L^2(\Omega)}^2.
\end{equation}
The definition of norm in $L^{\infty}(0,T;L^2(\Omega))$ gives:
\begin{equation*}
    \norm{q}_{L^{\infty}(0,T;L^2(\Omega))}^{2} \leq e^{\frac{16}{\mu} \norm{\mathbf{u}}^2_{\mathcal{U}}} \norm{q_0}_{L^2(\Omega)}^2
\end{equation*}
while squaring and integrating between $0$ and $T$ finally yields
\begin{equation*}
    \norm{q}_{L^{2}(0,T;L^2(\Omega))}^2 \leq T\,e^{\frac{16}{\mu} \norm{\mathbf{u}}_{\mathcal{U}}^2} \norm{q_0}_{L^2(\Omega)}^2.
\end{equation*}
\QEDA

\subsubsection{Proof of Estimate (\ref{subeq3})}
From Equation (\ref{int_norm}) we have:
\begin{equation*}
\begin{split}
        &\int_0^T \alpha_0(t) \norm{q(t)}^2_{H^1(\Omega)} \, dt  \\ &\leq \int_0^T \lambda(t) \norm{q(t)}^2_{L^2(\Omega)} dt -\frac{1}{2}\norm{q(T)}_{L^2(\Omega)}^2 + \frac{1}{2}\norm{q_0}_{L^2(\Omega)}^2
\end{split}
\end{equation*}
which implies
\begin{equation}
\label{build_est_3}
\int_0^T \alpha_0(t) \norm{q}^2_{H^1(\Omega)} \, dt  \leq \int_0^T \lambda(t) \norm{q}^2_{L^2(\Omega)} dt  + \frac{1}{2}\norm{q_0}_{L^2(\Omega)}^2.
\end{equation}

From Equation (\ref{base_est}) and the expression for $\lambda(t)$ in Equation (\ref{lambda_eq}) we have:
\begin{equation*}
    \lambda(t) \norm{q(t)}^2_{L^2(\Omega)} \leq \frac{8}{\mu} \sum_{i=1}^4 \norm{ u_i(t) }^2_{L^{\infty}(\Gamma_i)}  e^{\frac{16}{\mu} \norm{\mathbf{u}}^2_{\mathcal{U}}} \norm{q_0}_{L^2(\Omega)}^2;
\end{equation*}
then, integrating between $0$ and $T$:
\begin{equation*}
        \int_0^T \lambda(t) \norm{q(t)}^2_{L^2(\Omega)} \, dt \leq \frac{8}{\mu}\norm{\mathbf{u}}^2_{\mathcal{U}} e^{\frac{16}{\mu} \norm{\mathbf{u}}^2_{\mathcal{U}}} \norm{q_0}_{L^2(\Omega)}^2.
\end{equation*}
Plugging this result into Equation (\ref{build_est_3}) gives:
\begin{equation*}
\begin{split}
    & \int_0^T \alpha_0(t) \norm{q(t)}^2_{H^1(\Omega)} \, dt \\
    &\leq \Big(\frac{1}{2}+\frac{8}{\mu}\norm{\mathbf{u}}^2_{\mathcal{U}} e^{\frac{16}{\mu} \norm{\mathbf{u}}^2_{\mathcal{U}}}\Big)\norm{q_0}_{L^2(\Omega)}^2 .
\end{split}
\end{equation*}
Using the definition of $\bar{\alpha}_0$ in (\ref{bar_alpha_0}) we have:
\begin{equation*}
\begin{split}
&\bar{\alpha}_0\int_0^T \norm{q(t)}^2_{H^1(\Omega)} \, dt \leq
    \int_0^T \alpha_0(t) \norm{q(t)}^2_{H^1(\Omega)} \, dt \\ 
    &\leq \Big(\frac{1}{2}+\frac{8}{\mu}\norm{\mathbf{u}}^2_{\mathcal{U}} e^{\frac{16}{\mu} \norm{\mathbf{u}}^2_{\mathcal{U}}}\Big)\norm{q_0}_{L^2(\Omega)}^2 
\end{split}
\end{equation*}

and thus we obtain
\begin{equation*}
    \norm{q}^2_{L^2(0,T;H^1(\Omega))}  \leq \frac{1}{\bar{\alpha}_0}\Big(\frac{1}{2}+\frac{8}{\mu}\norm{\mathbf{u}}^2_{\mathcal{U}} e^{\frac{16}{\mu} \norm{\mathbf{u}}^2_{\mathcal{U}}}\Big)\norm{q_0}_{L^2(\Omega)}^2.
\end{equation*}
\QEDA

\subsubsection{Proof of Estimate (\ref{subeq4})} 
Using the duality pairing formalism we recast the state equation as
\begin{equation*}
    \langle \dot{q}(t),\phi\rangle_{H^{1}(\Omega)^*,H^1(\Omega)} + a(q(t),\phi) = 0 \qquad \forall \phi \in H^1(\Omega).
\end{equation*}
From the continuity of the bilinear form and Cauchy-Schwarz inequality we get
\begin{equation*}
\begin{split}
& |\langle \dot{q}(t),\phi\rangle_{H^{1}(\Omega)^*,H^1(\Omega)} | = |a(q(t),\phi)| \\ 
& \leq \mu \norm{\nabla q(t)}_{L^2(\Omega)}\norm{\nabla \phi}_{L^2(\Omega)} \\
&+ \norm{\mathbf{v}(t)}_{L^{\infty}(\Omega)^2} \norm{q(t)}_{L^2(\Omega)} \norm{\nabla \phi}_{L^2(\Omega)}  \\ 
& \leq \Big(\mu \norm{q(t)}_{H^1(\Omega)}+\norm{\mathbf{v}(t)}_{L^{\infty}(\Omega)^2} \norm{q(t)}_{L^2(\Omega)}\Big) \norm{\phi}_{H^1(\Omega)}.
\end{split}
\end{equation*}
Note that we have $\norm{\nabla \phi}_{L^{2}(\Omega)} \leq \norm{\phi}_{H^{1}(\Omega)}$.
From the definition of the dual space norm, we have:
\begin{equation*}
    \norm{\dot{q}(t)}_{H^{1}(\Omega)^{*}} \leq \Big(\mu \norm{q(t)}_{H^1(\Omega)}+\norm{\mathbf{v}(t)}_{L^{\infty}(\Omega)^2} \norm{q(t)}_{L^2(\Omega)}\Big). 
\end{equation*}
Squaring, integrating and using Cauchy inequality we have:
\begin{equation*}
\begin{split}
     &\int_0^T \norm{\dot{q}(t)}_{H^{1}(\Omega)^{*}}^2  \, dt  \\ 
     &\leq \int_0^T \Big(\mu \norm{q(t)}_{H^1(\Omega)}+\norm{\mathbf{v}(t)}_{L^{\infty}(\Omega)^2} \norm{q(t)}_{L^2(\Omega)}\Big)^2  \, dt \\ 
     &\leq   \int_0^T 2\mu^2 \norm{q(t)}_{H^1(\Omega)}^2 \, dt \\
     & + \int_0^T 2\norm{\mathbf{v}(t)}_{L^{\infty}(\Omega)^2}^2 \norm{q(t)}_{L^2(\Omega)}^2 \, dt \\
     & \leq 2 \mu^2 \norm{q}_{L^2(0,T;H^1(\Omega))}^2 \\
     &+ 2\norm{\mathbf{v}}_{L^2(0,T;L^{\infty}(\Omega)^2)}^2 \norm{q}_{L^{\infty}(0,T;L^2(\Omega))}^2 \\
     & \leq 2 \mu^2 \frac{1}{\bar{\alpha}_0}\Big(\frac{1}{2}+\frac{8}{\mu}\norm{\mathbf{u}}^2_{\mathcal{U}} e^{\frac{16}{\mu} \norm{\mathbf{u}}^2_{\mathcal{U}}}\Big)\norm{q_0}_{L^2(\Omega)}^2 \\
     &+ 16 \norm{\mathbf{u}}^2_{\mathcal{U}}
    e^{\frac{16}{\mu} \norm{\mathbf{u}}^2_{\mathcal{U}}} \norm{q_0}_{L^2(\Omega)}^2.      \\
\end{split}
\end{equation*}
Finally, we obtain
\begin{equation*}
\begin{split}
    &\norm{\dot{q}}_{L^2(0,T;H^{1}(\Omega)^{*})}^2  \\
    &\leq \Big(\frac{\mu^2}{\bar{\alpha}_0}+ 16\Big(1+\frac{\mu}{\bar{\alpha}_0}\Big)\norm{\mathbf{u}}^2_{\mathcal{U}} \, e^{\frac{16}{\mu} \norm{\mathbf{u}}^2_{\mathcal{U}}}\Big)\norm{q_0}_{L^2(\Omega)}^2.
\end{split}
\end{equation*}
\QEDA

\subsubsection{Proof of Theorem \ref{diff_c_t_s}}
Let us consider Equation (\ref{c_t_s_problem}), define the vector field $\bb{f} = \bb{v}_{\bb{h}}q$. Its weak form reads:
\begin{equation}
\label{weak_form_z}
\begin{split}
    &\int_{\Omega} \frac{\partial z(t)}{\partial t}\,\phi \, d\Omega \\
    &+ \int_{\Omega} \mu \nabla z(t) \cdot \nabla \phi - z(t)\,\bb{v}_{\bb{u}}(t) \cdot \nabla \phi \, d\Omega = \int_{\Omega} \mathbf{f}(t) \cdot \nabla \phi \, d\Omega
\end{split}
\end{equation}
that is,
\begin{equation*}
\begin{split}
    &\int_{\Omega} \frac{\partial z(t)}{\partial t}\,\phi d\Omega + a(q(t),\phi) \\
    &=  \int_{\Omega} \mathbf{f}(t) \cdot \nabla \phi \, d\Omega \quad \forall \phi \in H^{1}(\Omega), \,\, a.e. \,\, t \in (0,T). 
\end{split}
\end{equation*}
Subtituting $\phi = z(t)$ and using the weak coercivity of $a(q,\phi)$, Cauchy-Schwarz and Cauchy inequality, we obtain, for every $\epsilon>0$,

\begin{equation*}
\begin{split}
    &\int_{\Omega} \frac{\partial z(t)}{\partial t}\,z(t) d\Omega -\lambda(t) \norm{z(t)}_{L^2(\Omega)}^2\\  &
    +(\alpha_0(t)-\epsilon) \norm{z(t)}_{H^1(\Omega)}^2   \leq \frac{1}{4\epsilon}\norm{\mathbf{f}(t)}^2_{L^2(\Omega)^2}.
\end{split}
\end{equation*}
Choosing $\epsilon = \alpha_0(t)$ we have
\begin{equation}
\begin{split}
    \label{control_to_state_1}
     & \int_{\Omega} \frac{\partial z(t)}{\partial t}\,z(t) \, d\Omega  \leq \frac{1}{4 \, \alpha_0(t)}\norm{\mathbf{f}(t)}^2_{L^2(\Omega)^2} + \lambda(t) \norm{z(t)}_{L^2(\Omega)}^2 \\
     &\leq  \frac{1}{4 \, \bar{\alpha}_0}\norm{\mathbf{f}(t)}^2_{L^2(\Omega)^2} + \lambda(t) \norm{z(t)}_{L^2(\Omega)}^2.
\end{split}
\end{equation}
Equation (\ref{control_to_state_1}) can be rewritten as:
\begin{equation*}
    \frac{1}{2}\frac{d}{dt}\norm{z(t)}^2_{L^2(\Omega)} \leq  \frac{1}{4 \, \bar{\alpha}_0}\norm{\mathbf{f}(t)}^2_{L^2(\Omega)^2} + \lambda(t) \norm{z(t)}_{L^2(\Omega)}^2.
\end{equation*}

\noindent Integrating in time between $0$ and $t$, and considering that $z(0)=0$, we obtain:

\begin{equation*}
\begin{split}
    & \norm{z(t)}^2_{L^2(\Omega)} \\ 
    & \leq \int_0^t \frac{1}{2 \, \bar{\alpha}_0}\norm{\mathbf{f}(\tau)}^2_{L^2(\Omega)^2} \, d\tau +
    \int_0^t \, 2 \lambda(\tau) \norm{z(\tau)}_{L^2(\Omega)}^2 \, d\tau
\end{split}
\end{equation*}

\noindent then applying Gronwall's Lemma and the expression for $\lambda(t)$ in (\ref{lambda_eq}) we get:
\begin{equation}
\label{base_est_z}
    \begin{split}
    & \norm{z(t)}_{L^2(\Omega)}^2   \leq e^{\int_0^t 2\,\lambda(\tau) d\tau}\int_0^t \frac{1}{2 \bar{\alpha}_0} \norm{\mathbf{f}(\tau)}^2_{L^2(\Omega)^2}\,d\tau \, \\
    &\leq \frac{1}{2 \bar{\alpha}_0} e^{\frac{16}{\mu}    \norm{ \mathbf{u}}^2_{\mathcal{U}}} \norm{\mathbf{f}}^2_{L^2(0,T;L^2(\Omega)^2)}.
    \end{split}
\end{equation}
Integrating in time between $0$ and $T$ we finally obtain:
\begin{equation*}
\norm{z}_{L^2(0,T;L^2(\Omega))}^2 \leq \frac{T}{2 \bar{\alpha}_0} e^{\frac{16}{\mu}   \norm{ \mathbf{u}}^2_{\mathcal{U}}} \norm{\mathbf{f}}^2_{L^2(0,T;L^2(\Omega)^2)}.
\end{equation*}
We need to find a bound on $\norm{\mathbf{f}}^2_{L^2(0,T;L^2(\Omega))}$ based on the norms of state and control functions. Recalling the definition of $\bb{f} = \bb{v}_{\bb{h}} q$, we have:

\begin{equation*}
\begin{aligned}
    &\norm{\mathbf{f}}^2_{L^2(0,T;L^2(\Omega)^2)} = \int_0^T \norm{\bb{f}(t)}_{L^2(\Omega)^2}^2 dt \\
    &\leq \int_0^T \norm{\bb{v}_{\bb{h}}(t)}^2_{L^\infty(\Omega)^2}  \norm{q(t)}_{L^2(\Omega)}^2 dt \\
    &\leq \norm{\bb{v}_{\bb{h}}}^2_{L^2(0,T;L^\infty(\Omega)^2)}\norm{q}_{L^{\infty}(0,T;L^2(\Omega))}^2;
    \end{aligned}
\end{equation*}
then, using the results obtained in Lemma \ref{lemma_1} and Theorem \ref{theorem_1},
 we have:
\begin{equation*}
\norm{\mathbf{f}}^2_{L^2(0,T;L^2(\Omega))} \leq  8\norm{\mathbf{h}}^2_{\mathcal{U}}\,e^{\frac{16}{\mu} \norm{\mathbf{u}}^2_{\mathcal{U}}} \norm{q_0}_{L^2(\Omega)}^2
\end{equation*}
so that, finally, we obtain
\begin{equation*}
\begin{split}
    &\norm{z}_{L^2(0,T;L^2(\Omega))}^2 \\
    &\leq \frac{4\,T}{ \bar{\alpha}_0}  (e^{\frac{16}{\mu} \norm{\mathbf{u}}^2_{\mathcal{U}}})^2 \norm{q_0}_{L^2(\Omega)}^2 \norm{\mathbf{h}}^2_{\mathcal{U}} = C(\norm{\mathbf{u}}^2_{\mathcal{U}})\norm{\mathbf{h}}^2_{\mathcal{U}}.
\end{split}
\end{equation*}
Proceeding in a similar way as for the state equation, the following estimates hold:

\begin{subequations}
\begin{align*}
    &\norm{z}_{L^2(0,T;H^1(\Omega))} \\
    &\leq \frac{8}{\bar{\alpha}_0^2}\norm{\mathbf{h}}^2_{\mathcal{U}}\,e^{\frac{16}{\mu} \norm{\mathbf{u}}^2_{\mathcal{U}}} \norm{q_0}_{L^2(\Omega)}^2 \Big( 1 + \frac{16}{\mu }  e^{\frac{16}{\mu}    \norm{ \mathbf{u}}^2_{\mathcal{U}}}  \norm{\mathbf{u}}^2_{\mathcal{U}}  \Big),  \\
    &\norm{ \dot{z}}_{L^2(0,T;H^1(\Omega)^*)}^{2}   \\
    &\leq \frac{32}{\bar{\alpha}_0^2} \norm{\mathbf{h}}^2_{\mathcal{U}} e^{\frac{16}{\mu} \norm{\mathbf{u}}^2_{\mathcal{U}}} \norm{q_0}_{L^2(\Omega)}^2 \\
    &\Big( 16 \mu    e^{\frac{16}{\mu}    \norm{ \mathbf{u}}^2_{\mathcal{U}}}  \norm{\mathbf{u}}^2_{\mathcal{U}}  
+  4 \norm{\mathbf{u}}^2_{\mathcal{U}} \bar{\alpha}_0 +1+\mu^2  \Big). 
\end{align*}
\end{subequations}

In order to prove that the control-to-state map $q=\Xi(\ut)$ is Fréchet differentiable with directional derivative $z=\Xi'(\ut)\mathbf{h}$, we define $R = \Xi(\ut+\mathbf{h}) - \Xi(\ut) - z$ where  $\Xi(\ut+\mathbf{h})$ and $\Xi(\ut)$ are weak solutions of (\ref{weak_form}) with control actions $\ut+\mathbf{h}$ and $\ut$, respectively, whereas $z$ solves Equation (\ref{weak_form_z}). As a consequence, $R$ satisfies a.e. on $t \in (0,T)$:
\begin{equation}
\begin{split}
\label{weak_form_R}
&\int_{\Omega} \frac{\partial R(t)}{\partial t}\phi \, d\Omega \\
&+ \int_{\Omega} \mu \nabla R(t) \cdot \nabla \phi \, d\Omega - \int_{\Omega}  \mathbf{v}_{\ut+\mbf{h}}(t) R(t) \cdot \nabla \phi \, d\Omega  
\\  
&= \int_{\Omega} z(t)\,\mathbf{v}_{\mathbf{h}}(t)\cdot \nabla \phi \, d\Omega   \quad \forall \phi \in H^{1}(\Omega).
\end{split}
\end{equation}
Since Equation (\ref{weak_form_R}) has the same form of Equation (\ref{weak_form_z}), it holds that:
\begin{equation*}
\begin{split}
    &\norm{R}_{L^2(0,T;L^2(\Omega))}^2 \\ 
    &\leq \frac{T}{2 \bar{\alpha}_0} e^{{\frac{16}{\mu}}   \norm{ \mathbf{u} + \mathbf{h} }^2_{\mathcal{U}}} \norm{\mathbf{v}_{\bb{h}}}^2_{L^2(0,T;L^\infty(\Omega)^2)}\norm{z}_{L^{\infty}(0,T;L^2(\Omega))}^2,
\end{split}
\end{equation*}
while from Equation (\ref{base_est_z}) we obtain
\begin{equation*}
     \norm{z}_{L^{\infty}(0,T;L^2(\Omega))}^2
 \leq \frac{4}{ \bar{\alpha}_0} \Big(e^{\frac{16}{\mu}    \norm{ \mathbf{u}}^2_{\mathcal{U}}}\Big)^2 \norm{\mathbf{h}}^2_{\mathcal{U}} \norm{q_0}_{L^2(\Omega)}^2,
 \end{equation*}
 hence we have:
 \begin{equation*}
 \begin{split}
&\norm{R}_{L^2(0,T;L^2(\Omega))}^2 \\
&\leq \frac{16 \, T}{ \bar{\alpha}_0^2} e^{{\frac{16}{\mu}}   \norm{ \mathbf{u} + \mathbf{h} }^2_{\mathcal{U}}}  \Big(e^{\frac{16}{\mu}    \norm{ \mathbf{u}}^2_{\mathcal{U}}}\Big)^2 \norm{\mathbf{h}}^4_{\mathcal{U}} \norm{q_0}_{L^2(\Omega)}^2.
\end{split}
\end{equation*}
Proceeding in a similar way we obtain
\begin{equation*}
\begin{split}
&\norm{R}_{L^{2}\left(0, T ; H^{1}(\Omega)\right)}^{2} \\
&\leq \frac{32}{\bar{\alpha}_{0}^3}\norm{q_0}_{L^2(\Omega)}^2\Big(e^{\frac{16}{\mu}    \norm{ \mathbf{u}}^2_{\mathcal{U}}}\Big)^2\\
&\left(\frac{8 \norm{\bb{u}
+\bb{h}}^2_{\mathcal{U}}e^{\frac{16}{\mu}\norm{\bb{u}
+\bb{h}}^2_{\mathcal{U}}}}{\mu}+1\right)  \norm{\mathbf{h}}^4_{\mathcal{U}} 
\end{split}
\end{equation*}

\begin{equation*}
\begin{split}
& \lVert\dot{R}\rVert_{\Htd}^2 \\
&\leq  \frac{128}{ \bar{\alpha}_0}\norm{q_0}_{L^2(\Omega)}^2\Big(e^{\frac{16}{\mu}    \norm{ \mathbf{u}}^2_{\mathcal{U}}}\Big)^2\\
&\Big(\frac{2}{\bar{\alpha}_0}\Big(8 \norm{\bb{u}
+\bb{h}}^2_{\mathcal{U}} e^{\frac{16}{\mu} \norm{\bb{u}
+\bb{h}}^2_{\mathcal{U}}} \Big( \frac{2 \mu}{\alpha_0}+1 \Big)+ \frac{2\mu^2}{\alpha_0}\Big)+1\Big) \norm{\bb{h}}^4_{\mathcal{U}}
\end{split}
\end{equation*}
that can be compactly written as:
\begin{equation*}
\norm{R}_{L^{2}\left(0, T ; H^{1}(\Omega)\right)}^{2} \leq  C_1(\norm{\bb{u}}_{\mathcal{U}},\norm{\bb{h}}_{\mathcal{U}})\norm{\mathbf{h}}^4_{\mathcal{U}} 
\end{equation*}
and
\begin{equation*}
\lVert\dot{R}\rVert_{\Htd}^2 \leq  C_2(\norm{\bb{u}}_{\mathcal{U}},\norm{\bb{h}}_{\mathcal{U}})\norm{\mathbf{h}}^4_{\mathcal{U}} 
\end{equation*}

\noindent so that finally:

\begin{equation}
\begin{split}
\label{final_est_R}
&\norm{R}_{H^1(0,T;H^1(\Omega),H^1(\Omega)^*)}^2 \\
&= \norm{R}_{\Ht}^2 + \lVert\dot{R}\rVert_{\Htd}^2 \\
&\leq  C\norm{\bb{h}}_{\mathcal{U}}^4.
\end{split}
\end{equation}
where $C=\max\{C_1,C_2\}$ and has a finite value as $\norm{\mathbf{h}}_{\mathcal{U}} \to 0$.
As a consequence, from Equation (\ref{final_est_R}), we have that $\norm{R}_{H^1(0,T;H^1(\Omega),H^1(\Omega)^{*})} \to 0$ as $\norm{\mathbf{h}}_{\mathcal{U}} \to 0$, thus implying the Fréchet differentiability of the control-to-state map.
\QEDA

\section*{Acknowledgments}
The authors would like to thank Sandro Salsa (Politecnico di Milano) for his valuable comments and the anynonymous reviewers whose comments improved the paper.


\ifCLASSOPTIONcaptionsoff
  \newpage
\fi





\bibliographystyle{IEEEtran}
\bibliography{IEEEabrv,biblist_OR.bib}

\begin{thebibliography}{10}
\providecommand{\url}[1]{#1}
\csname url@rmstyle\endcsname
\providecommand{\newblock}{\relax}
\providecommand{\bibinfo}[2]{#2}
\providecommand\BIBentrySTDinterwordspacing{\spaceskip=0pt\relax}
\providecommand\BIBentryALTinterwordstretchfactor{4}
\providecommand\BIBentryALTinterwordspacing{\spaceskip=\fontdimen2\font plus
\BIBentryALTinterwordstretchfactor\fontdimen3\font minus
  \fontdimen4\font\relax}
\providecommand\BIBforeignlanguage[2]{{%
\expandafter\ifx\csname l@#1\endcsname\relax
\typeout{** WARNING: IEEEtran.bst: No hyphenation pattern has been}%
\typeout{** loaded for the language `#1'. Using the pattern for}%
\typeout{** the default language instead.}%
\else
\language=\csname l@#1\endcsname
\fi
#2}}
\renewcommand\BIBentryALTinterwordstretchfactor{4}

\bibitem{boundary_1}
S.~{Shahrokhi}, A.~{Mahadev}, and A.~T. {Becker}, ``Algorithms for shaping a
  particle swarm with a shared input by exploiting non-slip wall contacts,'' in
  \emph{2017 IEEE/RSJ International Conference on Intelligent Robots and
  Systems (IROS)}, 2017, pp. 4304--4311.

\bibitem{boundary_2}
S.~{Shahrokhi}, J.~{Shi}, B.~{Isichei}, and A.~T. {Becker}, ``Exploiting
  nonslip wall contacts to position two particles using the same control
  input,'' \emph{IEEE Transactions on Robotics}, vol.~35, no.~3, pp. 577--588,
  2019.

\bibitem{NIAC}
M.~Quadrelli, ``{NIAC} {P}hase {II} {O}rbiting {R}ainbows: {F}uture {S}pace
  {I}maging with {G}ranular {S}ystems,'' \emph{NASA Jet Propulsion Laboratory},
  2017.

\bibitem{DCD}
S.~Basinger, G.~Swartzlander, and M.~Quadrelli, ``Dynamics and {C}ontrol of a
  {D}isordered {S}ystem in {S}pace,'' \emph{AIAA SPACE 2013 Conference}, 2013.

\bibitem{granular_dyn}
M.~Quadrelli, S.~Basinger, and G.~Swartzlander, ``Multi-scale {D}ynamics,
  {C}ontrol, and {S}imulation of {G}ranular {S}pacecraft,'' \emph{European
  Community on Computational Methods in Applied Sciences (ECCOMAS) Thematic
  Conference, Zagreb, Croatia, July 1-4, , 2013}, 2013.

\bibitem{optics_gran}
S.~Basinger, D.~Palacios, M.~Quadrelli, and et~al., ``Optics of a granular
  imaging system (i.e. “orbiting rainbows”),'' \emph{Spie Optical
  Engineering + Applications}, 2015.

\bibitem{Unc_gran}
M.~Quadrelli, S.~Basinger, and E.~Sidick, ``Unconventional imaging with
  contained granular media,'' \emph{Proceedings Volume 10410, Unconventional
  and Indirect Imaging, Image Reconstruction, and Wavefront Sensing}, 2017.

\bibitem{Nature2020}
M.~Z. Miskin, A.~J. Cortese, K.~Dorsey, E.~P. Esposito, M.~F. Reynolds, Q.~Liu,
  M.~Cao, D.~A. Muller, P.~L. McEuen, and I.~Cohen, ``Electronically
  integrated, mass-manufactured, microscopic robots,'' \emph{Nature}, vol. 584,
  no. 7822, pp. 557--561, Aug 2020.

\bibitem{karthik_survey}
K.~Elamvazhuthi and S.~Berman, ``Mean-field models in swarm robotics: a
  survey,'' \emph{Bioinspiration {\&} Biomimetics}, vol.~15, no.~1, p. 015001,
  nov 2019.

\bibitem{Milutinovic2006ModelingAO}
D.~Milutinovic and P.~U. Lima, ``Modeling and optimal centralized control of a
  large-size robotic population,'' \emph{IEEE Transactions on Robotics},
  vol.~22, pp. 1280--1285, 2006.

\bibitem{gen_grad}
K.~{Rudd}, G.~{Foderaro}, and S.~{Ferrari}, ``A generalized reduced gradient
  method for the optimal control of multiscale dynamical systems,'' in
  \emph{52nd IEEE Conference on Decision and Control}, 2013, pp. 3857--3863.

\bibitem{nc_fod}
G.~{Foderaro} and S.~{Ferrari}, ``Necessary conditions for optimality for a
  distributed optimal control problem,'' in \emph{49th IEEE Conference on
  Decision and Control (CDC)}, 2010, pp. 4831--4838.

\bibitem{kart_1}
K.~{Elamvazhuthi} and S.~{Berman}, ``Optimal control of stochastic coverage
  strategies for robotic swarms,'' in \emph{2015 IEEE International Conference
  on Robotics and Automation (ICRA)}, 2015, pp. 1822--1829.

\bibitem{kart_2}
\BIBentryALTinterwordspacing
K.~Elamvazhuthi, H.~Kuiper, and S.~Berman, ``Pde-based optimization for
  stochastic mapping and coverage strategies using robotic ensembles,''
  \emph{Automatica}, vol.~95, pp. 356 -- 367, 2018. [Online]. Available:
  \url{http://www.sciencedirect.com/science/article/pii/S0005109818302991}
\BIBentrySTDinterwordspacing

\bibitem{karthik_bil}
K.~{Elamvazhuthi}, H.~{Kuiper}, M.~{Kawski}, and S.~{Berman}, ``Bilinear
  controllability of a class of advection–diffusion–reaction systems,''
  \emph{IEEE Transactions on Automatic Control}, vol.~64, no.~6, pp.
  2282--2297, 2019.

\bibitem{rob_dep}
T.~Zheng, Q.~Han, and H.~Lin, ``Deployment of robotic swarms via density
  feedback control,'' 2020.

\bibitem{lieberman}
G.~M. Lieberman, \emph{Second order parabolic differential equations}.\hskip
  1em plus 0.5em minus 0.4em\relax World Scientific, 1996.

\bibitem{Ito_97}
K.~Ito and K.~Kunisch, ``Estimation of the convection coefficient in elliptic
  equations,'' \emph{Inverse Problems}, vol.~13, no.~4, pp. 995--1013, aug
  1997.

\bibitem{joshi_05}
H.~R. Joshi, ``Optimal control of the convective velocity coefficient in a
  parabolic problem,'' \emph{Nonlinear Analysis: Theory, Methods \&
  Applications}, vol.~63, no.~5, pp. e1383 -- e1390, 2005, invited Talks from
  the Fourth World Congress of Nonlinear Analysts (WCNA 2004).

\bibitem{glowinski2021bilinear}
R.~Glowinski, Y.~Song, X.~Yuan, and H.~Yue, ``Bilinear optimal control of an
  advection-reaction-diffusion system,'' 2021.

\bibitem{Fleig2017}
A.~Fleig and R.~Guglielmi, ``{Optimal Control of the Fokker–Planck Equation
  with Space-Dependent Controls},'' \emph{Journal of Optimization Theory and
  Applications}, vol. 174, no.~2, pp. 408--427, 2017.

\bibitem{roy2018fokker}
S.~Roy, M.~Annunziato, A.~Borz{\`\i}, and C.~Klingenberg, ``A fokker--planck
  approach to control collective motion,'' \emph{Computational Optimization and
  Applications}, vol.~69, no.~2, pp. 423--459, 2018.

\bibitem{OT}
S.Bandyopadhyay and M.~Quadrelli, ``Optimal {T}ransport {B}ased {C}ontrol of
  {G}ranular {I}maging {S}ystem in {S}pace,'' \emph{International Workshop on
  Satellite Constellations and Formation Flying 17-17}, 2017.

\bibitem{chinese_ot}
K.~Zhang and Y.~Zhang, ``Optimal reconfiguration with collision avoidance for a
  granular spacecraft using laser pressure,'' \emph{Acta Astronautica}, vol.
  160, pp. 163 -- 174, 2019.

\bibitem{evans}
L.~C. Evans, \emph{{Partial Differential Equations}}.\hskip 1em plus 0.5em
  minus 0.4em\relax American Mathematical Society, 2010.

\bibitem{act_model}
M.~Sitti and D.~S. Wiersma, ``Pros and cons: Magnetic versus optical
  microrobots,'' \emph{Advanced materials (Deerfield Beach, Fla.)}, vol.~32,
  no.~20, p. e1906766, May 2020.

\bibitem{quarteroni_valli}
A.~Quarteroni and A.~Valli, \emph{Numerical Approximation of Partial
  Differential Equations}, 1st~ed.\hskip 1em plus 0.5em minus 0.4em\relax
  Springer Publishing Company, 2008.

\bibitem{brezis}
H.~Br{\'e}zis, \emph{Functional analysis, Sobolev spaces and partial
  differential equations}.\hskip 1em plus 0.5em minus 0.4em\relax Springer,
  2011, vol.~2, no.~3.

\bibitem{simon}
J.~Simon, ``Compact sets in the space ${L}^p({O},{T}; {B})$,'' \emph{Annali di
  Matematica Pura ed Applicata}, vol. 146, pp. 65--96, 01 1986.

\bibitem{fredi}
F.~Tr{\"o}ltzsch, \emph{Optimal Control of Partial Differential Equations:
  Theory, Methods, and Applications}, ser. Graduate Studies in
  Mathematics.\hskip 1em plus 0.5em minus 0.4em\relax American Mathematical
  Society, 2010.

\bibitem{Duprez2019}
\BIBentryALTinterwordspacing
M.~Duprez and P.~Lissy, ``{Bilinear local controllability to the trajectories
  of the Fokker-Planck equation with a localized control},'' 2021. [Online].
  Available: \url{http://arxiv.org/abs/1909.02831}
\BIBentrySTDinterwordspacing

\bibitem{noflux_barbu}
V.~Barbu, ``The controllability of fokker--planck equations with reflecting
  boundary conditions and controllers in diffusion term,'' \emph{SIAM Journal
  on Control and Optimization}, vol.~59, no.~1, pp. 709--726, 2021.

\bibitem{herzog}
R.~Herzog and K.~Kunisch, ``Algorithms for pde-constrained optimization,''
  \emph{GAMM-Mitteilungen}, vol.~33, no.~2, pp. 163--176, 2010.

\bibitem{schlog}
R.~Buchholz, H.~Engel, E.~Kammann, and F.~Tröltzsch, ``On the optimal control
  of the schlögl-model,'' \emph{Computational Optimization and Applications},
  vol.~56, 09 2013.

\bibitem{quart}
A.~Quarteroni, \emph{Numerical Models for Differential Problems}.\hskip 1em
  plus 0.5em minus 0.4em\relax Springer Publishing Company, 2016.

\bibitem{mexipopt}
E.~Bertolazzi, ``Matlab interface for ipopt,''
  \url{https://github.com/ebertolazzi/mexIPOPT/releases/tag/1.1.2)}, 2021.

\bibitem{ipopt}
A.~W{\"a}chter and L.~T. Biegler, ``On the implementation of an interior-point
  filter line-search algorithm for large-scale nonlinear programming,''
  \emph{Mathematical programming}, vol. 106, no.~1, pp. 25--57, 2006.

\bibitem{MQS}
A.~Manzoni, S.~Salsa, and A.~Quarteroni, \emph{Optimal Control of Partial
  Differential Equations}.\hskip 1em plus 0.5em minus 0.4em\relax Springer,
  2021, vol.~??

\end{thebibliography}
%

\begin{IEEEbiography}[{\includegraphics[width=1in,height=1.25in,clip,keepaspectratio]{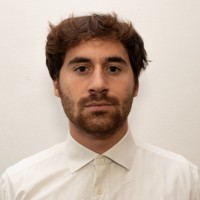}}]{Carlo Sinigaglia}
Carlo Sinigaglia is a Ph.D. Candidate in Mechanical Engineering at Politecnico di Milano. He obtained the M.Sc. in Mechanical Engineering from Politecnico di Milano, Italy, in 2019. His research interests include swarm robotics, macroscopic modeling of large-scale dynamical systems and optimal control.
\end{IEEEbiography}

\begin{IEEEbiography}[{\includegraphics[width=1in,height=1.25in,clip,keepaspectratio]{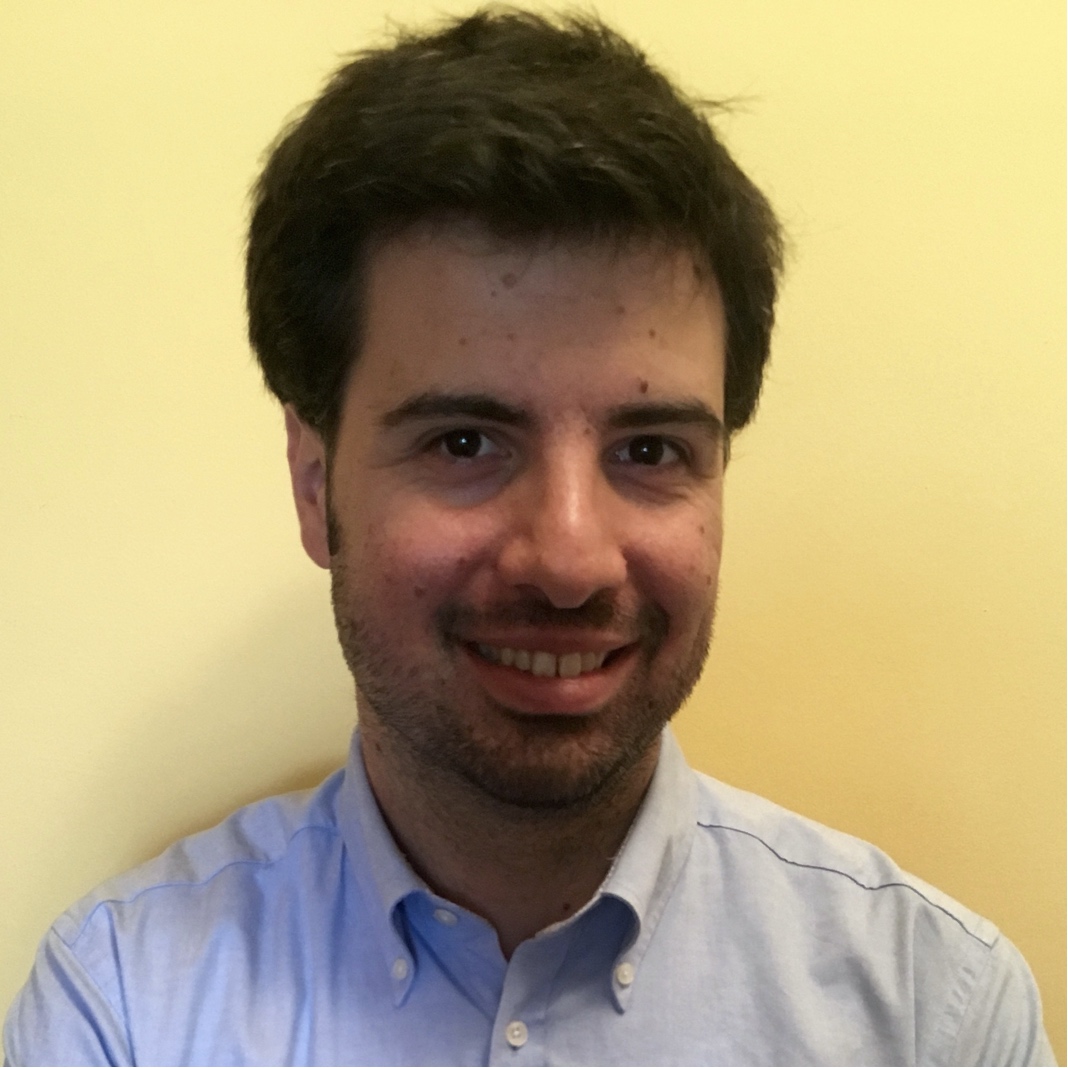}}]{Prof. Andrea Manzoni}
Andrea Manzoni received the M.Sc. Degree in Mathematical Engineering from Politecnico di Milano, Italy, in 2008, and the Ph.D. in Mathematics from the Ecole Polytechnique Fédérale de Lausanne (EPFL), Switzerland, in 2012. He has been research fellow at the International School of Advanced Studies, Trieste, Italy, until 2014, and researcher at EPFL until 2017, when he became tenure-track Assistant Professor at Politecnico di Milano. Since 2020 he is Associate Professor of Numerical Analysis at Politecnico di Milano.
\end{IEEEbiography}

\begin{IEEEbiography}[{\includegraphics[width=1in,height=1.25in,clip,keepaspectratio]{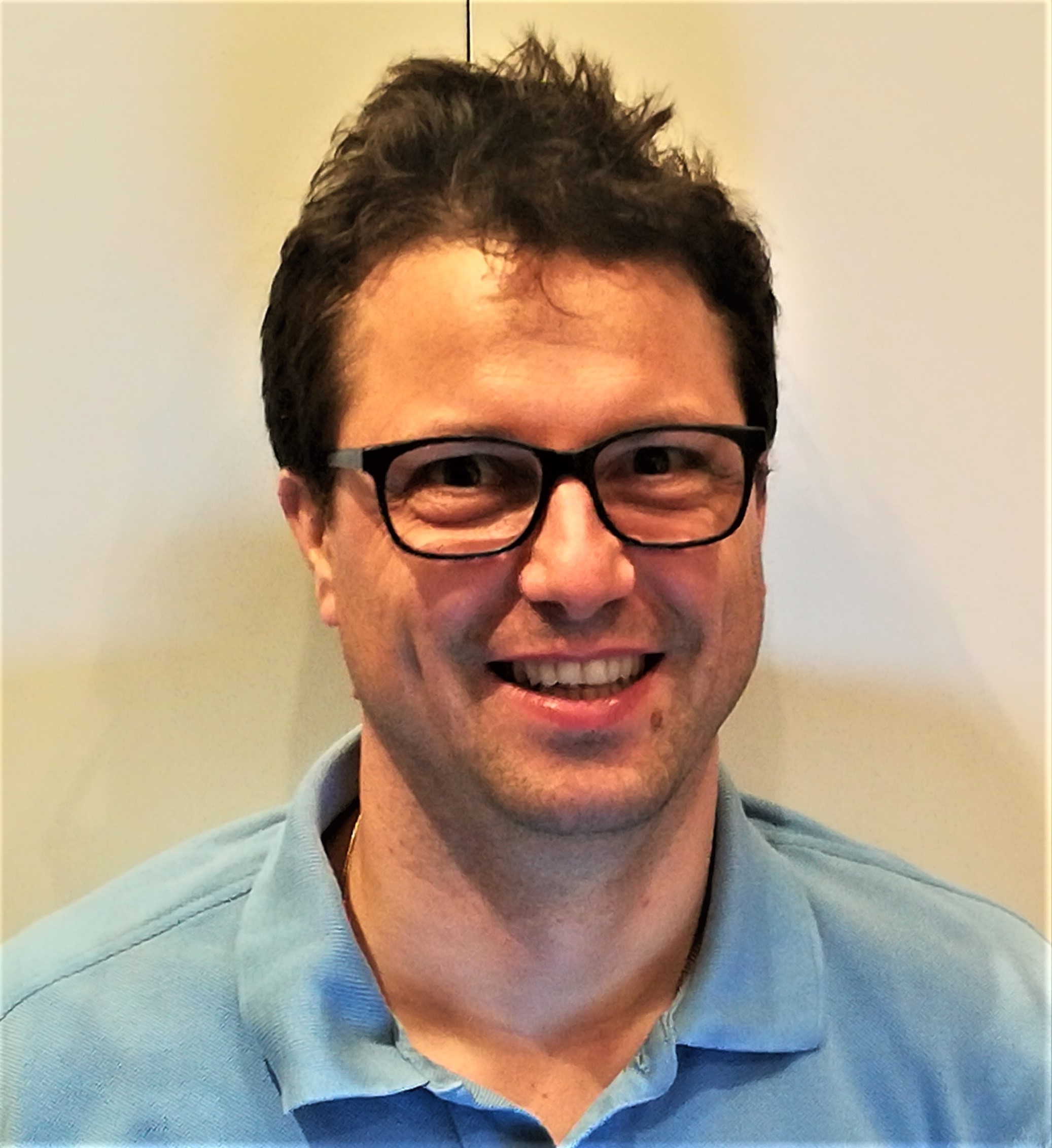}}]{Prof. Francesco Braghin}
Francesco Braghin received the M.Sc. degree in mechanical engineering and the Ph.D. degree in applied mechanics both from Politecnico di Milano, Italy, in 1997 and 2001, respectively. He became a Researcher, in 2001 and in 2011 an Associate Professor in the Department of Mechanical Engineering, Politecnico di Milano, where, since 2015, he has been a Full Professor in applied mechanics.
\end{IEEEbiography}

\end{document}